\documentclass{amsart}

\newenvironment{EllTab}{
\centerline
{Table 2. List of extremal elliptic $K3$ surfaces}
\par
\bigskip
\begin{tabular}{|c|l|l|ccc|} \hline 
 No & \hfil $\Sigma$ \hfil&\hfil  $MW$ \hfil& $a$ & $b$ & $c$ \\}
{\end{tabular}}

\newenvironment{SubConfigTable}{
\centerline
{Table 1. List of embedding of $\Delta$ in $\Gamma\sb f$}
\par
\bigskip
\begin{tabular}{|l||l|| c|l|c|} \hline 
\hfil no & \hfil $\Delta$ \hfil & No &\hfil  
$\Sigma\sb f$ \hfil & $\euler(\Sigma\sb f)$ \\}
{\end{tabular}}

\newtheorem{theorem}{Theorem}[section]
\newtheorem{lemma}[theorem]{Lemma}
\newtheorem{proposition}[theorem]{Proposition}
\newtheorem{corollary}[theorem]{Corollary}

\theoremstyle{definition}

\newtheorem{claim}{Claim}

\theoremstyle{remark}
\newtheorem{remark}[theorem]{\it Remark}
\newtheorem{step}{\it Step}

\numberwithin{equation}{section}

\renewcommand{\P}{\mathord{\mathbb P}}
\newcommand{\R}{\mathord{\mathbb R}}
\newcommand{\Z}{\mathord{\mathbb Z}}
\newcommand{\Q}{\mathord{\mathbb Q}}

\newcommand{\LLL}{\mathord{\mathcal L}}
\newcommand{\AAA}{\mathord{\mathcal A}}
\newcommand{\TTT}{\mathord{\mathcal T}}
\newcommand{\SSS}{\mathord{\mathcal S}}
\newcommand{\QQQ}{\mathord{\mathcal Q}}

\newcommand{\GL}{\mathord{GL}}
\newcommand{\SL}{\mathord{SL}}

\newcommand{\isomto}{
\hskip 2pt \smash{\mathop{\to}\limits\sp{\hskip -1pt \sim}}
\hskip 2pt}
\newcommand{\inv}{\sp{-1}}

\newcommand{\singtype}[1]{\mathord{\rm #1}}

\newcommand{\rt}{\operatorname{\rm root}\nolimits}
\newcommand{\rank}{\operatorname{\rm rank}\nolimits}
\newcommand{\euler}{\operatorname{\it eu\/}\nolimits}
\newcommand{\length}{\operatorname{\rm length}\nolimits}

\newcommand{\Aut}{\operatorname{\rm Aut}\nolimits}
\newcommand{\VERT}{\operatorname{\rm Vert}\nolimits}

\begin{document}

\title[extremal elliptic $K3$ surfaces]{Classification of extremal 
elliptic $K3$ surfaces
\\
and fundamental groups of open $K3$ surfaces}

\author{Ichiro Shimada}
\address{Department of Mathematics, Faculty of Science, Sapporo, JAPAN 061-0081}
\email{shimada@math.sci.hokudai.ac.jp}

\author{De-Qi Zhang}
\address{Department of Mathematics,  National University of Singapore,  
Lower KentRidge Road,
 SINGAPORE 119260}
\email{matzdq@math.nus.edu.sg}


\begin{abstract}
We present a complete list of extremal elliptic $K3$ surfaces (Theorem 1.1).
As an application,
we give a sufficient condition  for the topological fundamental group of
complement to an $ADE$-configuration of smooth rational curves on  a $K3$ surface
to be trivial (Proposition 4.1 and  Theorem 4.3).
\end{abstract}
\maketitle

\section{Introduction}
A complex elliptic $K3$ surface $f : X\to \P\sp 1$
with a section $O$
is said to be {\em extremal\/}
if the Picard number $\rho (X)$ 
of $X$ is $20$ and the Mordell-Weil group
$MW \sb f$ of $f$ is finite. The purpose of this paper is to
present the complete  list of all extremal elliptic $K3$
surfaces. 
As an application,
we show that, if
an $ADE$-configuration of smooth rational curves on a $K3$ surface
satisfies a certain condition,
then the topological fundamental group of the complement is  trivial.
(See Theorem~\ref{thm:2} for the precise statement.)
\par
\medskip
Let $f : X\to \P\sp 1$
be an elliptic $K3$ surface with a section $O$.
We denote by $R \sb f $ the set of all points $v \in \P\sp 1$
such that $f\sp{-1} (v)$ is reducible.
For a point $v\in R\sb f$,
let $f\inv( v) \sp\#$ be the union of irreducible components
of $f\sp{-1} (v)$
that are disjoint from the zero section  $O$.
It is known that the cohomology classes of irreducible 
components
of $f\inv( v) \sp\#$ form a negative definite root lattice $S\sb{f, v}$
of type $A\sb l$, $D\sb m$ or $E\sb n$
in $H\sp 2 (X; \Z)$.
Let $\tau (S\sb{f, v})$ be the type of this lattice.
We define $\Sigma \sb f$ to be the formal sum
of these types;
$$
\Sigma \sb f:=\sum\sb{v\in R\sb f} \tau (S\sb{f, v}).
$$ 
The N\'eron-Severi lattice  $NS \sb X$ of $X$
is defined to be $H\sp{1, 1} (X)\cap H\sp 2 (X; \Z)$,
and 
the transcendental lattice $T\sb X$
of $X$ is defined to be
the orthogonal complement of $NS \sb X$ in $H\sp 2 (X; \Z)$.
We call the triple $(\Sigma \sb f, MW \sb f, T\sb X)$
the {\em data} of the  elliptic $K3$ surface $f: X\to \P\sp 1$.
When $f : X\to \P\sp 1$ is extremal,
the transcendental lattice
$T\sb X$ is a positive definite even lattice of rank $2$.
\begin{theorem}
There exists
an extremal elliptic $K3$ surface  $f : X \to \P\sp 1$
with data $(\Sigma \sb f, MW \sb f, T \sb X)$
if and only if $(\Sigma \sb f, MW \sb f, T \sb X)$
appears in Table~{\rm 2} given at the end of this paper.
\end{theorem}
In Table~2,
the transcendental lattice $T\sb X$
is expressed by the coefficients of its Gram matrix
$$
\begin{pmatrix}
a & b \\ b & c 
\end{pmatrix}.
$$
See Subsection~\ref{subsec:2-1}
on how to  recover the $K3$ surface $X$ from $T\sb X$.
\par
\medskip
The classification of {\em semi-stable\/}
extremal elliptic $K3$ surfaces
has been done by
Miranda and Persson~\cite{MP}
and complemented by Artal-Bartolo,
Tokunaga and Zhang~\cite{ATZ}.
We can check that the semi-stable part 
of our list (No.~1- No.~112)
coincides with theirs.
Nishiyama~\cite{Nishiyama}
classified all elliptic fibrations
(not necessarily extremal)
on certain $K3$ surfaces.
On the other hand,
Ye~\cite{Ye} has independently
 classified all extremal elliptic $K3$ surfaces
with no semi-stable singular fibers
by different methods from ours.
\par
\medskip
\textbf{Acknowledgment.}
The authors would like to thank Professors 
Shigeyuki Kond\=o, 
Ken-ichi Nishiyama and Keiji Oguiso
for helpful discussions.
\section{Preliminaries}
\subsection{Transcendental lattice of singular $K3$ surfaces}\label{subsec:2-1}
Let $\QQQ$ be the set of symmetric  matrices
$$
Q=\begin{pmatrix}
a & b \\ b & c 
\end{pmatrix}
$$
of integer coefficients such that $a$ and $c$
are even and that the corresponding 
quadratic forms are positive definite.
The group $\GL \sb 2 (\Z)$ acts on $\QQQ$ from right by
$$
Q\mapsto {}\sp t\hskip -2pt g \cdot  Q \cdot g,
$$
where $g\in\GL\sb 2 (\Z)$.
Let $Q\sb 1$ and $Q\sb 2$ be two matrices in $\QQQ$,
and let $L\sb 1$ and $L\sb 2$
be the positive definite even lattices of rank  $2$ 
whose Gram matrices are $Q\sb1$ and $Q\sb 2$, respectively.
Then $L\sb 1$ and $L\sb 2$ are isomorphic as lattices
if and only if 
$Q\sb 1$ and $Q\sb 2$ are in the same orbit under the action of
$GL\sb 2 (\Z)$.
On the other hand,
each orbit in $\QQQ$ under the action of $\SL \sb 2 (\Z)$
contains a unique matrix with coefficients satisfying
$$
-a<2\;b\le a\le c, \qquad \textrm{with}\quad b\ge 0 \quad\textrm{if}\quad a=c.
$$
(See, for example, Conway and Sloane~\cite[p.~358]{CS}.)
Hence each orbit in $\QQQ$ under the action of $\GL \sb 2 (\Z)$
contains a unique matrix with coefficients satisfying
\begin{equation}\label{coef}
0\le 2\;b\le a\le c.
\end{equation}
In  Table~2, the transcendental lattice   is represented by the Gram matrix 
satisfying the condition~\eqref{coef}.
\par
Let $X$ be a $K3$ surface with $\rho (X)=20$;
that is, $X$ is a singular $K3$ surface in the terminology of Shioda and 
Inose~\cite{Shioda-Inose}.
The transcendental lattice $T\sb X$ can be naturally  oriented by means of
a holomorphic two form on $X$ (cf.~\cite[p.~128]{Shioda-Inose}).
Let $\SSS$ denote the set of  isomorphism classes of singular $K3$ surfaces.
Using the natural  orientation
on the transcendental lattice, we can lift the map $\SSS\to \QQQ/ \GL\sb 2 (\Z)$
given by $X\mapsto T\sb X$
to the map $\SSS \to \QQQ/\SL \sb 2 (\Z)$.
\begin{proposition}[Shioda and Inose~\cite{Shioda-Inose}]
This map $\SSS \to \QQQ/\SL \sb 2 (\Z)$
is bijective.
\qed
\end{proposition}
Moreover, Shioda and Inose~\cite{Shioda-Inose} gave us a method  to construct 
explicitly the singular $K3$ surface
corresponding to  a given element of $\QQQ/\SL \sb 2 (\Z)$
by means of Kummer surfaces.
The injectivity of the map $\SSS \to \QQQ/\SL \sb 2 (\Z)$
had been proved by Piateskii-Shapiro and  Shafarevich~\cite{PS}.
\par
\medskip
Suppose that an orbit $[ Q ] \in \QQQ/\GL\sb 2 (\Z)$
is represented by a matrix $Q$ satisfying~\eqref{coef}.
Let $\rho : \QQQ/\SL \sb 2 (\Z) \to \QQQ/\GL\sb 2 (\Z)$
be the natural projection.
Then we have
$$
|\rho\sp{-1} ([Q]) |=\begin{cases} 
2 & \textrm{if $0< 2\;b<a<c$} \\
1 & \textrm{otherwise.} 
\end{cases}
$$
Therefore, if a data in  Table~2 satisfies $a=c$ or $b=0$ or $2\;b=a$ 
(resp. $0<2\;b <a<c$),
then the number of the isomorphism classes of
  $K3$ surfaces that possess a structure of the extremal elliptic $K3$ surfaces  
with  the given data is one (resp. two). 
\subsection{Roots of  a negative definite even lattice}\label{subsec:roots}
Let $M$ be a negative definite even lattice.
A vector  of $M$ is said to be a {\em root} of $M$
if its norm is $-2$.
We denote by $\rt (M)$ the number of roots of $M$,
and by $M\sb{root}$ the sublattice of $M$
generated by the roots of $M$.
Suppose that a Gram matrix $ (a\sb{ij})$
of $M$ is given.
Then $\rt (M)$ can be  calculated by the following method.
Let
$$
g\sb r (x)=-\sum\sb{i, j=1}\sp r a\sb{ij}x\sb i x \sb j
$$
be the positive definite quadratic form associated with the opposite lattice $M\sp{-}$
of $M$,
where $r$ is the rank of $M$.
We consider the bounded closed subset
$$
E (g\sb r, 2):=\{ x \in \R\sp r \; ;\; g\sb r (x)\le 2 \}
$$
of $\R\sp r$.
Then we have
$$
\rt (M)+1=| E(g\sb r , 2) \cap\Z\sp r |,
$$
where $+1$ comes from the origin.
For a positive integer $k$ less than $r$, 
we write by  $p\sb k : \R\sp r \to \R\sp k$
 the projection $(x\sb 1, \dots, x\sb r)\mapsto (x\sb 1, \dots, x\sb k)$.
Then there exist a positive definite quadratic form $g\sb k$
of variables $(x\sb 1, \dots, x\sb k)$
and a positive  real number $\sigma\sb k$
such that 
$$
p\sb k ( E(g\sb r , 2))=E(g\sb k, \sigma\sb k):=
\{ y\in \R\sp k \; ; \; g\sb k (y)\le \sigma\sb k\}.
$$
The projection $(x\sb 1, \dots, x\sb{k+1})\mapsto (x\sb 1, \dots, x\sb k)$
maps $E(g\sb{k+1}, \sigma\sb{k+1})$ to $E(g\sb k, \sigma\sb k)$.
Hence,
if we have the list of the points of $E(g\sb k, \sigma\sb k)\cap\Z\sp{k}$,
then it is easy to make the list of the points of 
$E(g\sb{k+1}, \sigma\sb{k+1})\cap\Z\sp{k+1}$.
Thus,
starting from $E(g\sb 1, \sigma\sb 1)\cap \Z$,
we can make the list of the points of $E(g\sb r , 2) \cap\Z\sp r$
by induction on $k$.
\subsection{Root lattices of type $ADE$}
A {\em root type} is, by definition, 
a finite formal sum $\Sigma$ of $A\sb l$, $D\sb m$
and $E\sb n$ with non-negative integer coefficients;
$$
\Sigma =\sum\sb{l\ge 1} a\sb l A\sb l + 
\sum\sb{m\ge 4} d\sb m D\sb m
+\sum\sb {n=6}\sp 8 e\sb n E\sb n.
$$
We denote by $L (\Sigma)$ the negative definite root
lattice corresponding to $\Sigma$.
The rank of $L (\Sigma)$ is given by 
$$
\rank (L (\Sigma) ) = \sum\sb{l\ge 1} a\sb l  l + 
\sum\sb{m\ge 4} d\sb m  m
+\sum\sb {n=6}\sp 8 e\sb n  n,
$$
and the number of roots of $L (\Sigma)$ is given by
\begin{equation}\label{rootnumberADE}
\rt (L (\Sigma))=\sum\sb {l\ge 1} a\sb l  (l\sp 2 +l)+
\sum\sb{m\ge 4} d\sb m (2\,m\sp 2 - 2\, m)
+72\, e\sb 6 +126\, e\sb 7 + 240\, e\sb 8. 
\end{equation}
(See, for example, Bourbaki~\cite{Bourbaki}.)
Because of $L(\Sigma)\sb{root}=L(\Sigma)$,
we have
\begin{equation}\label{sigma}
L(\Sigma\sb 1)\cong L(\Sigma\sb 2) \Longleftrightarrow \Sigma\sb 1=\Sigma\sb 2.
\end{equation}
We also define  $\euler (\Sigma)$
by 
$$
\euler (\Sigma )
:= \sum\sb{l\ge 1} a\sb l  (l+1) + 
\sum\sb{m\ge 4} d\sb m  (m+2)
+\sum\sb {n=6}\sp 8 e\sb n  (n+2).
$$
\begin{lemma}\label{lem:euler}
Let $f : X\to\P\sp 1$ be an elliptic $K3$ surface.
Then $\euler (\Sigma\sb f)$ is at most $24$.
Moreover, if\/ $\euler (\Sigma\sb f)<24$,
then there exists  at least one singular fiber of type $\singtype{ I}\sb 1$,
$\singtype{ II}$, $\singtype{III}$ or $\singtype{IV}$.
\end{lemma}
\begin{proof}
Let $e(Y)$ denote the topological euler number
of a $CW$-complex $Y$.
Then $e(X)=24$ is equal with the sum of topological euler numbers of 
singular fibers of $f$.
Every singular fiber has a positive topological euler number.
We have defined $\euler (\Sigma)$ 
in such a way that, if $v \in R\sb f$,
then $\euler (\tau (S\sb{f, v}))\le e (f\sp{-1} (v)) $
holds,
and if $\euler (\tau (S\sb{f, v}))< e (f\sp{-1} (v)) $,
then the type of the fiber $f\sp{-1} (v)$ is either $\singtype{III}$ or $\singtype{IV}$.
Hence $\euler (\Sigma\sb f)$ does not exceed the sum of the topological
euler numbers of reducible singular fibers,
and if $\euler (\Sigma\sb f)<24$, then there is  an irreducible singular fiber
or a singular fiber of type $\singtype{III}$ or $\singtype{IV}$.
\end{proof}
\subsection{Discriminant form and overlattices}
Let $L$ be an even lattice,
$L\sp\vee$ the dual of $L$,
$D\sb L $ the discriminant group $L\sp\vee /L$
of $L$,
and $q\sb L$ the discriminant form on $D\sb L$.
(See Nikulin~\cite[n.~4]{Nikulin2}
for the definitions.)
An overlattice of $L$ is, by definition, an integral 
sublattice of the $\Q$-lattice $L\sp{\vee}$
containing $L$.
\begin{lemma}[Nikulin~\cite{Nikulin2} Proposition 1.4.2]\label{lem:Nikulin1}
{\rm (1)}
Let $A$ be an isotopic subgroup of $(D\sb L, q\sb L)$.
Then the pre-image $M:=\phi\sb L \sp{-1} (A)$
of $A$ by the natural projection $\phi\sb L : L\sp{\vee}\to D\sb L$
is an overlattice of $L$,
and the discriminant form $(D\sb M, q\sb M)$ of $M$
is isomorphic to 
$(A\sp{\perp}/A, q\sb L |\sb{A\sp{\perp}/A})$,
where $A\sp{\perp}$ is the orthogonal complement
of $A$ in $D\sb L$,
and $q\sb L |\sb{A\sp{\perp}/A}$ is
the restriction of $q\sb L$ to $A\sp{\perp}/A$.
{\rm (2)}
The correspondence $A\mapsto M$
gives a bijection 
from the set of isotopic subgroups 
of $(D\sb L, q\sb L)$ to
the set of even overlattices of $L$.
\qed
\end{lemma}
\begin{lemma}[Nikulin~\cite{Nikulin2} Corollary 1.6.2]\label{lem:Nikulin2}
Let $S$ and $K$ be two even lattices.
Then the following two conditions are equivalent.
{\rm (i)}
There is an isomorphism $\gamma : D\sb S \isomto D\sb K$
of abelian groups
such that $\gamma\sp * q\sb K =-q\sb S$.
{\rm (ii)}
There is an even unimodular overlattice 
of $S\oplus K$ into which $S$ and $K$ are primitively
embedded.
\qed
\end{lemma}
\subsection{N\'eron-Severi groups of elliptic $K3$ surfaces}
Let $f : X\to \P\sp 1$ be an elliptic $K3$ surface
with the zero section $O$.
In the N\'eron-Severi lattice  $NS\sb X$ of $X$,
the cohomology classes of the zero section 
 $O$ and a general fiber of $f$ generate a sublattice $U \sb f$
of rank $2$,
which is isomorphic to the hyperbolic lattice
$$
H:=\begin{pmatrix} 0 & 1 \\ 1 & 0 \end{pmatrix}.
$$
Let $W\sb f$ be the orthogonal complement of $U\sb f$ in $NS\sb X$.
Because $U\sb f$ is unimodular,
we have $NS\sb X = U\sb f \oplus W\sb f$.
Because $U\sb f$ is of signature $(1,1)$
and $NS\sb X$ is of signature $(1, \rho (X)-1)$,
$W\sb f$ is negative definite of rank $\rho (X)-2$.
Note that $W\sb f$ contains the sublattice
$$
S\sb f:=\bigoplus \sb {v\in R\sb f} S\sb{f, v}
$$
generated by the cohomology classes of irreducible components
of reducible fibers of $f$ that are disjoint from the zero section.
By definition,
$S\sb f$ is isomorphic to $L (\Sigma \sb f)$.
\begin{lemma}[Nishiyama~\cite{Nishiyama} Lemma 6.1]\label{lem:Nishiyama}
The sublattice $S\sb f$ of $W\sb f$ coincides with $(W\sb f)\sb{root}$,
and the Mordell-Weil group $MW\sb f$ of $f$
is isomorphic to $W\sb f / S\sb f$.
In particular,
$\rt (L (\Sigma\sb f))$ is equal with $\rt (W\sb f)$.
\qed
\end{lemma}
Because $W\sb f \oplus U\sb f \oplus T\sb X$
has an even unimodular overlattice $H\sp 2 (X; \Z)$
into which $NS\sb X=W\sb f \oplus U\sb f$
and $T\sb X$
are primitively embedded,
and because the discriminant form of $NS\sb X$
is equal with the discriminant form of $W\sb f$
by $D\sb{ U\sb f}=(0)$,
Lemma~\ref{lem:Nikulin2} implies the following:
\begin{corollary}\label{cor:q}
There is an isomorphism
$\gamma : D \sb{W\sb f}\isomto  D\sb {T\sb X}$
of abelian groups
such that $\gamma\sp * q\sb{T\sb X}$
coincides with $-q\sb{W\sb f}$. 
\qed
\end{corollary}
\subsection{Existence of elliptic $K3$ surfaces}
Let $\Lambda$ be the $K3$ lattice $L (2E\sb 8) \oplus H\sp{\oplus 3}$.
\begin{lemma}[Kond\=o~\cite{Kondo1} Lemma 2.1]\label{lem:Kondo}
Let $T$ be a positive definite primitive sublattice
of $\Lambda$
with $\rank (T)= 2$,
and $T\sp\perp$ the orthogonal complement of $T$ in $\Lambda$.
Suppose that $T\sp{\perp}$
contains a sublattice $H\sb T$
isomorphic to the hyperbolic lattice.
Let $M\sb T$ be the orthogonal complement of $H\sb T$ in $T\sp{\perp}$.
Then there exist an elliptic $K3$ surface
$f : X\to \P\sp 1$
such that $T\sb X \cong T$
and $W\sb f \cong M\sb T$.
\end{lemma}
\begin{proof}
By the surjectivity of the period map
of the moduli of $K3$ surfaces (cf. Todorov~\cite{Todorov}),
there exist a $K3$ surface $X$ and an isomorphism
$\alpha : H\sp 2 (X; \Z)\cong \Lambda$ of lattices
such that $\alpha\sp{-1} (T)=T\sb X$.
By  Kond\=o~\cite[Lemma 2.1]{Kondo1},
the $K3$ surface $X$ has an elliptic fibration $f : X\to \P\sp 1$
with a section such that $\Z [F]\sp{\perp} /\Z [F] \cong M\sb T$,
where $[F]\in U\sb f$ is the cohomology class of a fiber of $f$,
and $\Z [F]\sp{\perp}$ is the orthogonal complement of $[F]$ in the 
N\'eron-Severi lattice 
$NS\sb X$.
Because  $NS\sb X $ is equal with $U\sb f \oplus W\sb f$,
and because $\Z [F]\sp{\perp} \cap U\sb f$ coincides with $\Z [F]$,
we see that $\Z [F]\sp{\perp} / \Z [F] $ is isomorphic to $W\sb f$.
\end{proof}
\subsection{Datum of extremal elliptic $K3$ surfaces}
\begin{proposition}\label{prop:data}
A triple $(\Sigma, MW, T)$
consisting of a root type  $\Sigma$,
a finite abelian group $MW$
and a positive definite even lattice $T$ of rank $2$
is a data of 
an extremal elliptic $K3$ surface 
if and only if the following hold\/{\rm :}
\begin{enumerate}
\renewcommand{\labelenumi}{$(D\arabic{enumi})$}
\item
$\length (MW)\le 2$,  $\rank ( L(\Sigma))=18$ and $\euler (\Sigma)\le 24$.
\item
 There exists an overlattice $M$ of $L(\Sigma)$
satisfying the following\/{\rm :}
\begin{enumerate}
\renewcommand{\labelenumii}{$(D2\:\text{-}\:\alph{enumii})$}
\item $M /L (\Sigma)\cong MW$,
\item there exists  an isomorphism $\gamma : D\sb{M}\isomto D\sb T$
of abelian groups
such that $\gamma\sp * q\sb T=-q\sb{M}$,
and
\item $\rt (L (\Sigma))=\rt (M)$.
\end{enumerate}
\end{enumerate}
\end{proposition}
\begin{proof}
Suppose that there exists an extremal elliptic $K3$ surface
$f : X\to \P\sp 1$
with data equal with $(\Sigma, MW, T)$.
It is obvious that $\Sigma$ and $MW$ satisfies  the condition $(D1)$.
Via the isomorphism $S\sb f \cong L (\Sigma)$,
the overlattice $W\sb f$ of $S\sb f$
corresponds to an overlattice $M$ of $L (\Sigma)$,
which satisfies the conditions $(D2\:\text{-}\:a)$-$(D2\:\text{-}\:c)$
by Lemma~\ref{lem:Nishiyama} and Corollary~\ref{cor:q}.
Conversely, suppose that $(\Sigma, MW, T)$
satisfies the conditions $(D1)$ and $(D2)$.
By Lemma~\ref{lem:Nikulin2},
the condition $(D2\:\text{-}\:b)$ and $D\sb H =0$ imply that there exists an even
unimodular  overlattice of $M \oplus  H \oplus T$
into which $M \oplus H$ and $T$ are
primitively embedded.
By the theorem of Milnor (see,  for example, Serre~\cite{Serre})
on the classification of even unimodular lattices,
any even unimodular lattice of signature $(3, 19)$
is isomorphic to the $K3$ lattice $\Lambda$.
Then Lemma~\ref{lem:Kondo} implies that 
there exists an elliptic $K3$ surface 
$f : X\to \P\sp 1$ satisfying $W\sb f \cong M$
and $T\sb X\cong T$.
The condition $(D2\:\text{-}\:c)$
implies $M \sb{root}=L(\Sigma)$.
Combining this with Lemma~\ref{lem:Nishiyama},
we see that $S\sb f \cong L (\Sigma)$.
Then (2.2) implies that $\Sigma\sb f=\Sigma$.
Using Lemma~\ref{lem:Nishiyama} and the condition $(D2\:\text{-}\:a)$,
we see that $MW\sb f \cong MW$.
Thus the data of $f : X\to \P\sp 1$ coincides with $(\Sigma, MW, T)$.
\end{proof}
\begin{remark}
In the light of Lemma~\ref{lem:Nikulin1},
the condition $(D2)$ 
is equivalent to the following:
\begin{enumerate}
\setcounter{enumi}{2}
\renewcommand{\labelenumi}{$(D\arabic{enumi})$}
\item
There exists an isotopic subgroup
$A$ of $(D\sb{L(\Sigma)}, q\sb{L(\Sigma)})$
satisfying the following:
\begin{enumerate}
\renewcommand{\labelenumii}{$(D3\:\text{-}\:\alph{enumii})$}
\item 
$A$ is isomorphic to $MW$,
\item 
there exists an isomorphism 
$\gamma : A\sp{\perp}/ A\isomto D\sb T$ of abelian groups 
such that $\gamma\sp * q\sb T =-q\sb{L(\Sigma)} |\sb{ A\sp{\perp}/A}$, and 
\item
$\rt (\phi\sb{L(\Sigma)}\sp{-1} (A))$ is equal with $\rt (L(\Sigma))$,
where $\phi\sb{L(\Sigma)} : L(\Sigma )\sp{\vee} \to D\sb{L(\Sigma)}$
is the natural projection.
\end{enumerate}
\end{enumerate}
\end{remark}
\begin{remark}
We did not use the conditions $\length (MW)\le 2$ and
$\euler (\Sigma)\le 24$ in the proof of 
the ``\thinspace if\thinspace" part of Proposition~\ref{prop:data}.
It follows that,
if $(\Sigma, MW, T)$ satisfies $\rank (L (\Sigma))=18$
and the condition $(D2)$,
then $\length (MW)\le 2$ and $\euler (\Sigma)\le 24$
follow automatically.
This fact can be used 
when we check the computer program described in the next section.
\end{remark}
\section{Making the list}
First we list up
all root types $\Sigma$
satisfying $\rank (L (\Sigma))=18$ and $\euler (\Sigma)\le 24$.
This list $\LLL$ consists of $712$ elements.
\par
\medskip
Next
we run a program that takes an element
$\Sigma$ of the list $\LLL$
as an input and proceeds as follows.
\begin{step}
The program calculates the intersection matrix of $L(\Sigma)\sp{\vee}$.
Using this matrix,
it calculates the discriminant form of $L(\Sigma)$,
and decomposes it into $p$-parts;
$$
(D\sb{L(\Sigma)}, q\sb{L(\Sigma)})=
\bigoplus \sb{p} (D\sb{L(\Sigma)}, q\sb{L(\Sigma)})\sb p,
$$
where $p$ runs through the set $\{p\sb 1, \dots, p\sb k\}$
of prime divisors of the discriminant
$| D\sb{L(\Sigma)}|$ of $L(\Sigma)$.
We write the $p\sb i$-part of $(D\sb{L(\Sigma)}, q\sb{L(\Sigma)})$
by $(D\sb{L(\Sigma), i}, q\sb{L(\Sigma), i})$.
\end{step}
\begin{step}
For each $p\sb i$,
it calculates the set
$I( p\sb i)$ of all pairs $(A, A\sp{\perp})$
of an isotopic subgroup $A$
of $(D\sb{L(\Sigma), i}, q\sb{L(\Sigma), i})$
and its orthogonal complement $A\sp{\perp}$
such that $\length (A)\le 2$.
\end{step}
\begin{step}
For each element
$$
\AAA:=((A\sb 1, A\sb 1\sp{\perp}), \dots, (A\sb k, A\sb k\sp{\perp}) )
\in I (p\sb 1)\times\cdots\times I (p\sb k),
$$
it calculates the $\Q/2\Z$-valued quadratic  form
$$
q\sb{\AAA}:=
q\sb{L(\Sigma), 1} |\sb{A\sb 1\sp{\perp}/A\sb 1}
\times\cdots\times
q\sb{L(\Sigma), k} |\sb{A\sb k\sp{\perp}/A\sb k}
$$
on the finite abelian group
$$
D\sb{\AAA}:=A\sb 1\sp{\perp}/A\sb 1
\times\cdots\times A\sb k\sp{\perp}/A\sb k.
$$
Let $d(\AAA)$ be the order of $D\sb {\AAA}$.
\end{step}
\begin{step}
It generates the list $\TTT (d(\AAA))$
of positive definite even lattices of rank $2$
with discriminant equal with $d(\AAA)$.
For each $T\in \TTT (d(\AAA))$,
it calculates the discriminant form of $T$
and decomposes it into $p$-parts.
If $D\sb T$ is isomorphic to $D\sb {\AAA}$
and $q\sb T$ is isomorphic to $-q\sb{\AAA}$,
then it proceeds to the next step.
Note that the automorphism group
of a finite abelian $p$-group of length $\le 2$
is easily calculated,
and hence it is an easy task
to check whether
two given quadratic forms
on the finite abelian $p$-group of length $\le 2$ are isomorphic  or not.
\end{step}
\begin{step}
It calculates the Gram matrix
of the sublattice $\widetilde L (\AAA)$
of $L (\Sigma)\sp{\vee}$ 
generated by $L (\Sigma)\subset L (\Sigma)\sp{\vee}$ and the
pull-backs of generators of the subgroups 
$A\sb i\subset D\sb {L(\Sigma),  i}$
by the projection 
$L(\Sigma)\sp{\vee} \to D\sb {L(\Sigma)}\to D\sb {L(\Sigma), i}$.
Then it calculates $\rt (\widetilde L (\AAA))$
by the method described in the subsection~\ref{subsec:roots}.
If $\rt (\widetilde L (\AAA))$
is equal with $\rt (L(\Sigma))$
calculated by~\eqref{rootnumberADE},
then it puts out the pair of the finite abelian group
$$
MW:=A\sb 1 \times \cdots \times A\sb k
$$
and the lattice $T$.
\end{step}
Then $(\Sigma, MW, T)$ satisfies the conditions $(D1)$ and $(D3)$,
and all triples $(\Sigma, MW, T)$ satisfying $(D1)$ and $(D3)$
are obtained by this program.
\section{Fundamental groups of open $K3$ surfaces}
A simple normal crossing divisor $\Delta$ on a $K3$ surface $X$
is said to be an {\em $ADE$-configuration of smooth rational curves}    
if each irreducible component of $\Delta$ is a smooth
rational curve
and the intersection matrix of the irreducible
components
of $\Delta$ is a direct sum
of the Cartan matrices of type $A\sb l$, $D\sb m$
or $E\sb n$ multiplied by $-1$.
It is known that $\Delta$ is an $ADE$-configuration of smooth rational curves   
if and only if each connected component of
$\Delta$ can be contracted to a rational double point.
We consider the following quite plausible hypothesis.
Let $\Delta$ be an $ADE$-configuration of smooth rational curves
on a $K3$ surface $X$.
\par
\medskip
\noindent
{\bf Hypothesis.}
If $\pi\sb 1\sp{alg} (X\setminus  \Delta)$ is trivial,
then so is $\pi\sb 1 (X\setminus  \Delta)$.
\par
\medskip\noindent
Here $\pi\sb 1 \sp{alg} (X\setminus  \Delta)$ is the algebraic fundamental 
group of $X\setminus  \Delta$,
which is the pro-finite completion of
the topological fundamental group $\pi\sb 1 (X\setminus  \Delta)$.
\begin{proposition}\label{prop:from-hyp}
Suppose that Hypothesis is true
for any $ADE$-configuration of smooth rational curves   
on an arbitrary $K3$ surface.
Let $\Delta$ be an $ADE$-configuration of smooth rational curves   
on a $K3$ surface $X$.
Then $\pi\sb 1 (X\setminus  \Delta)$
satisfies one of the following\/{\rm :}
\begin{enumerate}
\renewcommand{\labelenumi}{(\roman{enumi})}
\item
$\pi\sb 1 (X\setminus  \Delta)$ is trivial.
\item
There exist a complex torus $T$ of dimension $2$ 
and a finite automorphism group $G$ of $T$
such that $T/G$ is birational to $X$
and that $\pi\sb 1 (X\setminus  \Delta)$
fits in the exact sequence
$$
1\longrightarrow \pi\sb 1 (T) \longrightarrow 
\pi\sb 1 (X\setminus  \Delta) \longrightarrow G \longrightarrow 1.
$$
\item
$\pi\sb 1 (X\setminus  \Delta)$ is isomorphic to a 
symplectic automorphism group of a $K3$ surface.
\end{enumerate}
\end{proposition}
\begin{remark}
Fujiki~\cite{Fujiki} classified the automorphism groups of
complex tori of dimension $2$.  In particular,
the  $G$  in (ii) is either one of  $\Z/(n)$ ($n = 2, 3, 4, 6$),
$Q_8$ (Quaternion of order 8), $D_{12}$ (Dihedral of order 12)
and  $T_{24}$ (Tetrahedral of order 24), whence
the  $\pi_1(X \setminus \Delta)$  in (ii) is a soluble group.
Mukai~\cite{Mukai} presented the complete list
of symplectic automorphism groups of
$K3$ surfaces. (See also Kond\=o~\cite{Kondo2} and Xiao~\cite{Xiao}.)
Under Hypothesis, therefore,
we know what groups can appear as $\pi\sb 1 (X\setminus  \Delta)$.
\end{remark}
\textit{Proof of Proposition~{\rm \ref{prop:from-hyp}}.}
Suppose that $\pi\sb 1 (X\setminus  \Delta)$ is non-trivial.
By Hypothesis,
$\pi\sb 1 \sp{alg} (X\setminus  \Delta)$ is also
non-trivial.
For a surjective homomorphism $\phi : \pi\sb 1 (X\setminus \Delta)\to G$
from $\pi\sb 1 (X\setminus \Delta)$ to a finite group $G$,
we denote by
$$
\psi\sb{\phi} : \widetilde Y\sb{\phi} \longrightarrow X
$$
the finite Galois cover of $X$ corresponding to $\phi$,
which is \'etale over $X\setminus \Delta$
and whose Galois group is canonically isomorphic to $G$.
Let $\rho : \widetilde Y\sb{\phi}\sp\prime \to \widetilde Y\sb{\phi}$ be the 
resolution of singularities,
and  $\gamma : \widetilde Y\sb{\phi}\sp\prime \to Y\sb\phi$ 
the contraction of $(-1)$-curves.
We denote by $\Delta\sb\phi$ the union of 
one-dimensional irreducible components of 
$\gamma(\rho\sp{-1} (\psi\sb{\phi}\sp{-1}(\Delta)))$.
Then it is easy to see that $Y\sb\phi$ is either a $K3$ surface or a complex 
torus of dimension $2$,
and that the Galois group $G$ of $\psi\sb\phi$ acts on $Y\sb\phi$ symplectically.
Moreover, 
$\Delta\sb\phi$ is an empty set or an $ADE$-configuration of
smooth rational curves.
We have an exact sequence
$$
1
\;\longrightarrow\; \pi\sb 1 (Y\sb\phi\setminus\Delta\sb \phi)
\;\longrightarrow\; \pi\sb 1 (X\setminus\Delta)
\;\longrightarrow\; G
\;\longrightarrow\; 1,
$$
because $\pi\sb 1 (\widetilde Y\sb{\phi}\setminus \psi\sb{\phi}\sp{-1} (\Delta))$
is isomorphic to $\pi\sb 1 (Y\sb{\phi}\setminus\Delta\sb{\phi})$.
Suppose that there exists $\phi : \pi\sb 1 (X\setminus \Delta)\to G$
such that $Y\sb\phi$ is a complex torus 
of dimension $2$.
Then $\Delta\sb \phi$ is empty, and hence (ii) occurs.
Suppose that no complex tori of dimension $2$
appear as a finite Galois cover of $X$ branched in $\Delta$.
Then any finite quotient group of $\pi\sb 1 (X\setminus  \Delta)$
must appear in Mukai's list of symplectic 
automorphism groups of $K3$ surfaces.
Because this list consists 
of finite number of isomorphism classes
of finite groups,
there exists a maximal finite quotient
$\phi\sb{max} : \pi\sb 1 (X\setminus  \Delta)\to G\sb{max}$
of $\pi\sb 1 (X\setminus  \Delta)$.
Then $\pi\sb 1 (Y\sb{\phi\sb{max}}\setminus \Delta\sb{\phi\sb{max}})$
has no non-trivial finite quotient group,
and hence
it is trivial by Hypothesis.
Thus (iii) occurs.
\qed
\par
\medskip
For an $ADE$-configuration $\Delta$ of smooth rational curves   
on a $K3$ surface $X$,
we denote by $\Z [\Delta]$ the sublattice of $H\sp 2 (X; \Z)$
generated by the cohomology classes of
the irreducible components of $\Delta$,
which is isomorphic to a negative definite root lattice of type $ADE$.
We denote by $\Sigma\sb{\Delta}$ the  root type
such that $\Z [\Delta]$ is isomorphic to $L (\Sigma\sb{\Delta})$.
Using the list of extremal elliptic $K3$ surfaces,
we prove the following theorem.
We first consider the following conditions on
a root type $\Sigma$ (see (2.4) for the definition of  $D\sb {L(\Sigma)}$).
\par
\medskip
\begin{enumerate}
\renewcommand{\labelenumi}{$(N\arabic{enumi})$}
\item
$\rank (L (\Sigma))\le 18$, and 
\item
$\length ( D\sb {L(\Sigma)})\le 20-\rank (L(\Sigma))$.
\end{enumerate}
\begin{theorem}\label{thm:2}
Let  $X$  be a  $K3$  surface and  $\Delta$  an  $ADE$-configuration
of smooth rational curves on  $X$.  
Suppose that the root type $\Sigma\sb\Delta$ satisfies 
conditions $(N1)$ and $(N2)$. 
If $\Z [\Delta]$ is primitive in $H\sp 2 (X; \Z)$
then $\pi\sb 1 (X\setminus \Delta)$ is trivial.
\end{theorem}
In virtue of Lemma~\ref{lem:Xiao} below,
we can easily derive the following:
\begin{corollary}
Let  $X$  be a  $K3$  surface and  $\Delta$  an  $ADE$-configuration
of smooth rational curves on  $X$.  
Suppose that  $\Sigma\sb{\Delta}$  satisfies the conditions 
$(N1)$  and  $(N2)$.
Then Hypothesis is true.
\qed
\end{corollary}
\begin{remark}
The conditions $(N1)$ and $(N2)$ come from Nikulin~\cite[Theorem~1.14.1]{Nikulin2}
(see also Morrison~\cite[Theorem~2.8]{Morrison}),
which gives a sufficient condition for the uniqueness of 
the primitive embedding of $L (\Sigma)$ into the $K3$ lattice
$\Lambda$.
\end{remark}
First we prepare some  lemmas.
Let $\overline{\Z[\Delta]}$ be the primitive closure 
of $\Z[\Delta]$ in $H\sp 2 (X; \Z)$.
\begin{lemma}[Xiao~\cite{Xiao} Lemma 2]\label{lem:Xiao}
The dual of 
the abelianisation of $\pi\sb 1 (X\setminus  \Delta)$ is canonically isomorphic to
$\overline{\Z[\Delta]}/\Z[\Delta]$. In particular,
if $\pi\sb 1 \sp{alg} (X\setminus  \Delta)$ is trivial,
then $\Z [\Delta]$ is primitive in $H\sp 2 (X; \Z)$.
\qed
\end{lemma}
Let $\Gamma\sb 1$ and $\Gamma\sb 2$ be graphs
with the set of vertices denoted by
$\VERT (\Gamma\sb 1)$ and  $\VERT(\Gamma\sb2)$,
respectively.
An embedding of $\Gamma\sb 1$ into $\Gamma\sb 2$
is, by definition,
an injection $f : \VERT (\Gamma\sb 1)\to \VERT(\Gamma\sb2)$
such that, for any $u, v \in \VERT (\Gamma\sb 1)$,
$f(u)$ and $f(v)$ are connected by an edge of $\Gamma\sb 2$
if and only if $u$ and $v$ are connected by an edge of $\Gamma\sb 1$.
\par
\medskip
Let $\Gamma (\Sigma)$ denote the Dynkin graph of $\Sigma$.
\begin{lemma}\label{lem:ext}
Suppose that $\Sigma$ satisfies the conditions $(N1)$ and  $(N2)$.
Then there exists $\Sigma\sp\prime$
satisfying $\rank (L (\Sigma\sp\prime))=18$ and the condition  $(N2)$
such that $\Gamma (\Sigma)$ can be  embedded in $\Gamma (\Sigma\sp\prime)$.
\end{lemma}
\begin{proof}
This  is  checked by listing up all $\Sigma$
satisfying the conditions $(N1)$ and $(N2)$ using computer.
\end{proof}
\begin{lemma}\label{lem:IIIandIV}
Let $f : X\to \P\sp 1$ be an elliptic surface with the zero section $O$.
Suppose that a fiber $f\sp{-1} (v)$ over $v\in \P\sp 1$ is a singular fiber of type
$\singtype{III}$ or $\singtype{IV}$.
Let $\Xi$ be a union of some irreducible components
of $f\sp{-1} (v)$ that does not coincide with the whole fiber $f\sp{-1} (v)$.
If $U$ is a small open disk on $\P\sp 1$ with the center $v$,
then $f\sp{-1} (U)\setminus (\Xi\cup (f\sp{-1} (U) \cap O))$
has an abelian fundamental group.
\end{lemma}
\begin{proof}
This can be proved easily by the van-Kampen theorem.
\end{proof}
\begin{lemma}\label{lem:NM}
Let $\Sigma$ be satisfying  the conditions  $(N1)$ and  $(N2)$.
Suppose that  $(X, \Delta)$ and $(X\sp\prime, \Delta\sp\prime)$
satisfy the following{\rm:}
\begin{enumerate}
\renewcommand{\labelenumi}{(\alph{enumi})}
\item $\Sigma\sb{\Delta}=\Sigma\sb{\Delta\sp\prime}=\Sigma$,
\item $\overline{\Z[\Delta]}=\Z [\Delta]$ and $\overline{\Z[\Delta\sp\prime]}=\Z [\Delta\sp\prime]$.
\end{enumerate}
Then there exists a connected continuous family
$(X\sb t, \Delta\sb t)$ parameterized by $t\in [0, 1]$
such that $(X\sb 0, \Delta\sb 0)=(X, \Delta)$,
$(X\sb 1, \Delta\sb 1)=(X\sp\prime, \Delta\sp\prime)$
and that $(X\sb t, \Delta\sb t)$ are diffeomorphic to one another.
In particular, $\pi\sb 1 (X\setminus  \Delta)$ is isomorphic
to $\pi\sb 1 (X\sp\prime\setminus  \Delta\sp\prime)$.
\end{lemma}
\begin{proof}
By Nikulin~\cite[Theorem~1.14.1]{Nikulin2},
the primitive embedding of $L (\Sigma)$ into the $K3$ lattice
$\Lambda$ is unique up to $\Aut (\Lambda)$.
Hence the assertion follows from Nikulin's
connectedness theorem~\cite[Theorem~2.10]{Nikulin1}.
\end{proof}
\textit{Proof of Theorem~{\rm \ref{thm:2}}.}
Let us consider the following:
\begin{claim}\label{claim:1}
Suppose that $\Sigma$ satisfies $\rank (L (\Sigma))=18$
and the condition $(N2)$.
Then there exists an $ADE$-configuration of smooth rational curves   
$\Delta \sb \Sigma$ on a $K3$ surface $X\sb \Sigma$
such that $\Sigma\sb{\Delta\sb \Sigma}=\Sigma$ and 
$\pi\sb 1 (X\sb \Sigma \setminus  \Delta\sb \Sigma)=\{ 1 \}$.
\end{claim}
We deduce Theorem~\ref{thm:2} from Claim 1.
Suppose that $\Delta$ is an $ADE$-configuration of smooth rational 
curves   on a $K3$ surface $X$
such that $\Sigma\sb{\Delta}$ satisfies 
the conditions $(N1)$ and $(N2)$, 
and that 
$\Z [\Delta]$ is primitive in $H\sp 2 (X;\Z)$.
By Lemma~\ref{lem:ext},
there exists $\Sigma\sb 1$
satisfying $\rank (L(\Sigma\sb 1))=18$ and the condition $(N2)$
such that $\Gamma (\Sigma\sb{\Delta})$ is embedded into 
$\Gamma (\Sigma\sb 1)$.
By Claim 1, we have $(X\sb 1, \Delta \sb 1)$
such that  $\Sigma\sb{\Delta\sb 1}=\Sigma\sb 1$ and 
$\pi\sb 1 (X\sb 1\setminus  \Delta\sb 1)=\{ 1 \}$.
Let $\Delta\sp\prime\subset\Delta\sb 1 $ be
the sub-configuration of smooth rational curves on $X\sb 1$   
corresponding to the
subgraph $\Gamma (\Sigma\sb{\Delta})\hookrightarrow
\Gamma (\Sigma\sb 1)=\Gamma(\Sigma\sb{\Delta\sb 1})$.
There is a surjection from
$\pi\sb 1 (X\sb 1\setminus  \Delta\sb 1)$
to
$\pi\sb 1 (X\sb 1\setminus  \Delta\sp\prime)$,
and hence $\pi\sb 1 (X\sb 1\setminus  \Delta\sp\prime)$
is trivial.
In particular, $\Z [\Delta \sp\prime]$ is
primitive in $H\sp 2 (X\sb 1; \Z)$.
Because of $\Sigma\sb{\Delta\sp\prime}=\Sigma\sb{\Delta}$,
Lemma~\ref{lem:NM}
implies that $\pi\sb 1 (X\setminus  \Delta)$ is
isomorphic to $\pi\sb 1 (X\sb 1 \setminus  \Delta\sp\prime)$.
Thus $\pi\sb 1 (X\setminus  \Delta)$ is trivial.
\par
\medskip
Let $f : X\to \P\sp 1$ be an extremal elliptic $K3$ surface.
For a point $v\in R\sb f$,
we denote the total fiber of $f$ over $v$ by
$$
\sum\sb{i=1}\sp{r\sb v} m\sb{v, i} \: C\sb{v, i},
$$
where $m\sb{v, i}$ is the multiplicity of the irreducible 
component $C\sb{v, i}$ of $f\sp{-1} (v)$.
We denote by $\Gamma\sb f$ the union of the zero section
and all irreducible fibers $f\sp{-1} (v)$ $(v\in R\sb f)$.
\begin{claim}
Suppose that $MW\sb f=(0)$.
Suppose that a sub-configuration $\Delta$ of $\Gamma\sb f$
satisfies the following two conditions.
\begin{enumerate}
\renewcommand{\labelenumi}{$(Z\arabic{enumi})$}
\item
The number of $v\in R\sb f$ such that 
$m\sb{v, i}=1 \Longrightarrow C\sb{v, i} \subset \Delta$
holds  is at most one.
\item
Either one of the following holds:
\begin{enumerate}
\renewcommand{\labelenumii}{$(Z2\:\text{-}\:\alph{enumii})$}
\item
The configuration $\Delta$ does not contain the zero section, 
\item
there is a point $v\sb 1 \in R\sb f$
such that the type $\tau (S\sb{f, v\sb 1})$ is $A\sb 1$
and that $F\sb 1:=f\sp{-1} (v\sb 1)$ and $\Delta$ have no common
irreducible components, or
\item
$\euler (\Sigma\sb f) \le 23$.
\end{enumerate}
\end{enumerate}
Then $\pi\sb 1 (X\setminus \Delta)$ is trivial.
\end{claim}
\textit{Proof of Claim~2.}
By Lemma~\ref{lem:Nishiyama},
the assumption $MW\sb f=(0)$
implies that the cohomology classes $[O]$ and 
$[C\sb{v, i}]$ $ (v\in R\sb f, i=1, \dots, r\sb v)$ of the irreducible
components  of $\Gamma\sb f$ span $NS \sb X$.
The relations among these generators are generated by
$$
\sum\sb{i=1}\sp{r\sb v} m\sb{v, i} \: C\sb{v, i} =
\sum\sb{i=1}\sp{r\sb {v\sp{\prime}}} m\sb{v\sp{\prime}, i} \: C\sb{v\sp{\prime}, i}
\quad (v, v\sp{\prime}\in R\sb f).
$$
Therefore the  condition $(Z1)$
implies that the cohomology classes of
the irreducible components of $\Delta$
constitute a subset of a $\Z$-basis of $NS\sb X$.
Hence $\Z[\Delta]$ is primitive in $H\sp 2 (X; \Z)$.
In particular, $\pi\sb 1 (X\setminus \Delta)$ is a perfect group
by Lemma~\ref{lem:Xiao}.
On the other hand,
the condition $(Z1)$ implies
that there exists a point $v\sb 0\in \P\sp 1$
such that every fiber of the restriction
$$
f|\sb{X\setminus (\Delta\cup f\sp{-1} (v\sb 0)) } : 
X\setminus (\Delta\cup f\sp{-1} (v\sb 0)) \longrightarrow \P\sp 1 \setminus \{ v\sb 0 \}
$$
of $f$
has a reduced irreducible component.
Then, by Nori's lemma~\cite[Lemma 1.5 (C)]{Nori},
 if $U$ is a non-empty connected classically open subset of 
$\P\sp 1\setminus \{v\sb 0 \}$,
then the inclusion of 
$f\sp{-1} (U)\setminus (f\sp{-1} (U)\cap \Delta))$ into 
$X\setminus (\Delta \cup f\sp{-1} (v\sb 0))$
induces a surjection on the fundamental groups.
The inclusion of $X\setminus (\Delta \cup f\sp{-1} (v\sb 0))$
into $X\setminus\Delta$
also induces a surjection on the fundamental groups.
We shall show that
there exists a small open disk $U$ on $\P\sp 1 \setminus \{v\sb 0 \}$
such that
$$
G\sb U :=\pi\sb 1 (f\sp{-1} (U)\setminus (f\sp{-1} (U)\cap \Delta))
$$
is abelian.
When $(Z2\:\text{-}\:a)$ occurs, we take a small open disk disjoint from $R\sb f$ as
$U$. Then $G\sb U$
is  abelian, because of $f\sp{-1} (U)\cap \Delta=\emptyset$.
Suppose that $(Z2\:\text{-}\:b)$ occurs.
We can take $v\sb 0$ from $\P\sp 1 \setminus \{v\sb 1\}$,
because $F\sb 1$ has no irreducible components of 
multiplicity $\ge 2$.
We choose as $U$ a small open disk with the center $v\sb 1$.
There is a contraction from $f\sp{-1} (U)\setminus (f\sp{-1} (U)\cap\Delta)$
to $F\sb 1 \setminus (F\sb 1 \cap \Delta)$.
Because $\pi\sb 1(F\sb 1 \setminus (F\sb 1 \cap \Delta))$ is abelian, so is $G\sb U$.
Suppose that $(Z2\:\text{-}\:c)$ occurs.
By Lemma~\ref{lem:euler},
there exists a singular fiber $F\sb 2 :=f\sp{-1} (v\sb 2)$
of type $\singtype{I}\sb 1$, $\singtype{II}$, $\singtype{III}$ or $\singtype{IV}$.
Because $F\sb 2$ has no irreducible components of multiplicity $\ge 2$,
we can choose $v\sb 0$ from $\P\sp 1 \setminus \{v\sb 2\}$.
If $F\sb 2 $ is of type $\singtype{I}\sb 1$ or $\singtype{II}$,
then $F\sb 2\cap \Delta$ consists of a nonsingular point of $F\sb 2$,
and $\pi\sb 1 (F\sb 2 \setminus (F\sb 2\cap \Delta))$ is abelian.
Hence $G\sb U$ is also abelian.
If $F\sb 2$ is of type $\singtype{III}$ or $\singtype{IV}$,
then $F\sb 2 \cap \Delta$ cannot coincide with the whole fiber $F\sb 2$.
Hence Lemma~\ref{lem:IIIandIV} implies that $G\sb U$ is abelian.
Therefore we see that   $\pi\sb 1 (X\setminus \Delta)$ is abelian.
Being both perfect and abelian,  $\pi\sb 1(X\setminus\Delta)$ is trivial.
\qed
\par
\medskip
Now we proceed to the proof of Claim~\ref{claim:1}.
We list up all
$\Sigma$ 
satisfying the condition $(N2)$
and $\rank (L (\Sigma))=18$.
It consists of $297$ elements.
Among them, $199$ elements 
can be the type $\Sigma\sb f$
of singular fibers 
of some extremal elliptic $K3$ surface
$f : X\to \P\sp 1$
with $MW\sb f=0$.
For these configurations,
$\pi\sb 1 (X\setminus \Delta)$ is trivial by Claim~2.
The remaining $98$ configurations are listed in 
the second column of Table~1 below.
Each of them is a sub-configuration of $\Gamma\sb f$
satisfying the conditions $(Z1)$ and $(Z2)$,
where $f : X \to \P\sp 1$ is the 
extremal elliptic $K3$ surface
with $MW\sb f=0$
whose number in Table~2
is given  in the third column of Table~1.
The fourth and fifth  columns of Table~1 indicate $\Sigma\sb f$
and $\euler (\Sigma\sb f)$, respectively.
In the case nos.~20, 28, 39, 41 and  85 in Table 1,
we can choose the embedding of $\Delta$ into $\Gamma\sb f$
in such a way that $(Z2\:\text{-}\:b)$ holds.
In the case nos.~30, 37, 57 and 63 in Table 1,
we can choose the embedding of $\Delta$ into $\Gamma\sb f$
in such a way that $(Z2\:\text{-}\:a)$ holds.
By Claim 2 again,
$\pi\sb 1 (X\setminus \Delta)$ is  trivial for 
these $98$ configurations $\Delta$.
\qed
\par
\medskip
\begin{remark}
The graph $\Gamma (A\sb{19})$ (resp.~$\Gamma (D\sb{19})$)
can be embedded into $\Gamma\sb f$
in such a way that $(Z1)$ and $(Z2)$ are satisfied,
where $f : X\to \P\sp 1$ is the extremal elliptic $K3$ surfaces
whose number in Table~2 is 312 (resp.~320).
Therefore, if $\Gamma (\Delta)$ is embedded in $\Gamma (A\sb{19})$
or $\Gamma (D\sb{19})$, then $\Gamma (\Delta)$
can be embedded in  $\Gamma\sb f$
in such a way that $(Z1)$ and $(Z2)$ are satisfied.
\end{remark}
\vfill\eject
\begin{SubConfigTable} 
\hline 1&$A\sb{2}+A\sb{3}+2\;A\sb{4}+A\sb{5}$&19&$A\sb{2}+2\;A\sb{3}+A\sb{4}+A\sb{6}$&$23$\\ 
\hline 2&$A\sb{1}+A\sb{2}+A\sb{3}+2\;A\sb{6}$&23&$A\sb{1}+A\sb{2}+A\sb{4}+A\sb{5}+A\sb{6}$&$23$\\ 
\hline 3&$2\;A\sb{1}+A\sb{4}+2\;A\sb{6}$&23&$A\sb{1}+A\sb{2}+A\sb{4}+A\sb{5}+A\sb{6}$&$23$\\ 
\hline 4&$2\;A\sb{2}+2\;A\sb{4}+A\sb{6}$&23 & $A\sb{1}+A\sb{2}+A\sb{4}+A\sb{5}+A\sb{6}$&$23$\\ 
\hline 5&$A\sb{1}+A\sb{5}+2\;A\sb{6}$&40&$A\sb{1}+A\sb{4}+A\sb{6}+A\sb{7}$&$22 $\\ 
\hline 6&$A\sb{4}+2\;A\sb{7}$&52&$A\sb{4}+A\sb{6}+A\sb{8}$&$21$\\ 
\hline 7&$A\sb{1}+A\sb{2}+2\;A\sb{4}+A\sb{7}$&23&$A\sb{1}+A\sb{2}+A\sb{4}+A\sb{5}+A\sb{6}$&$23$ \\ 
\hline 8&$A\sb{3}+2\;A\sb{4}+A\sb{7}$&24&$A\sb{3}+A\sb{4}+A\sb{5}+A\sb{6}$&$22$\\ 
\hline 9&$A\sb{2}+2\;A\sb{4}+A\sb{8}$&36&$A\sb{2}+A\sb{4}+A\sb{5}+A\sb{7}$&$22$\\ 
\hline 10&$2\;A\sb{3}+A\sb{4}+A\sb{8}$&46&$A\sb{1}+A\sb{2}+A\sb{3}+A\sb{4}+A\sb{8}$&$23$\\ 
\hline 11&$A\sb{3}+A\sb{7}+A\sb{8}$&53&$A\sb{1}+A\sb{2}+A\sb{7}+A\sb{8}$&$22$\\ 
\hline 12&$A\sb{1}+2\;A\sb{2}+A\sb{4}+A\sb{9}$&46 & $A\sb{1}+A\sb{2}+A\sb{3}+A\sb{4}+A\sb{8}$&$23$\\ 
\hline 13&$A\sb{2}+A\sb{3}+A\sb{4}+A\sb{9}$&71&$2\;A\sb{2}+A\sb{4}+A\sb{10}$&$22$\\ 
\hline 14&$A\sb{3}+A\sb{4}+A\sb{11}$&93&$A\sb{2}+A\sb{4}+A\sb{12}$&$21 $\\ 
\hline 15&$A\sb{7}+A\sb{11}$&312&$A\sb{10}+E\sb{8}$&$21$\\ 
\hline 16&$2\;A\sb{3}+A\sb{12}$&93 &$ A\sb{2}+A\sb{4}+A\sb{12}$&$21$\\ 
\hline 17&$A\sb{3}+A\sb{15}$&312 &$ A\sb{10}+E\sb{8}$&$21$\\ 
\hline 18&$A\sb{2}+2\;A\sb{6}+D\sb{4}$&99&$A\sb{2}+A\sb{3}+A\sb{13}$&$21$\\ 
\hline 19&$2\;A\sb{4}+A\sb{6}+D\sb{4}$&18&$A\sb{1}+A\sb{3}+2\;A\sb{4}+A\sb{6}$&$23$\\ 
\hline 20&$2\;A\sb{2}+A\sb{4}+A\sb{6}+D\sb{4}$&20&$A\sb{1}+2\;A\sb{2}+A\sb{3}+A\sb{4}+A\sb{6}$&$24$\\ 
\hline 21&$A\sb{2}+A\sb{4}+A\sb{8}+D\sb{4}$&44&$2\;A\sb{1}+2\;A\sb{4}+A\sb{8}$&$23$\\ 
\hline 22&$A\sb{6}+A\sb{8}+D\sb{4}$&50&$2\;A\sb{1}+A\sb{2}+A\sb{6}+A\sb{8}$&$23$\\ 
\hline 23&$2\;A\sb{2}+A\sb{10}+D\sb{4}$&72&$2\;A\sb{1}+A\sb{2}+A\sb{4}+A\sb{10}$&$23$\\ 
\hline 24&$A\sb{4}+A\sb{10}+D\sb{4}$&72&$2\;A\sb{1}+A\sb{2}+A\sb{4}+A\sb{10}$&$23$\\ 
\hline 25&$A\sb{2}+A\sb{12}+D\sb{4}$&90&$2\;A\sb{1}+2\;A\sb{2}+A\sb{12}$&$23$\\ 
\hline 26&$A\sb{14}+D\sb{4}$& 320&$D\sb{10}+E\sb{8}$&$22$\\ 
\hline 27&$2\;A\sb{2}+A\sb{4}+2\;D\sb{5}$&210&$2\;A\sb{2}+D\sb{14}$&$22$\\ 
\hline 28&$A\sb{1}+2\;A\sb{2}+2\;A\sb{4}+D\sb{5}$&157&$A\sb{1}+A\sb{2}+2\;A\sb{4}+D\sb{7}$&$24$\\ 
\hline 29&$A\sb{2}+A\sb{3}+2\;A\sb{4}+D\sb{5}$&46&$A\sb{1}+A\sb{2}+A\sb{3}+A\sb{4}+A\sb{8}$&$23$\\ 
\hline 30 &$ A\sb{2}+A\sb{6}+2\;D\sb{5}$&193&$A\sb{2}+A\sb{6}+D\sb{10}$&$22$\\ 
\hline 31&$A\sb{3}+A\sb{4}+A\sb{6}+D\sb{5}$&18&$A\sb{1}+A\sb{3}+2\;A\sb{4}+A\sb{6}$&$23$\\ 
\hline 32&$A\sb{2}+A\sb{4}+A\sb{7}+D\sb{5}$&72&$2\;A\sb{1}+A\sb{2}+A\sb{4}+A\sb{10}$&$23$\\ 
\hline 33&$A\sb{6}+A\sb{7}+D\sb{5}$&50&$2\;A\sb{1}+A\sb{2}+A\sb{6}+A\sb{8}$&$23$\\ 
\hline 34&$A\sb{2}+A\sb{3}+A\sb{8}+D\sb{5}$&50&$2\;A\sb{1}+A\sb{2}+A\sb{6}+A\sb{8}$&$ 23$\\ 
\hline 35&$A\sb{3}+A\sb{10}+D\sb{5}$&69&$A\sb{1}+2\;A\sb{2}+A\sb{3}+A\sb{10}$&$23$\\ 
\hline 36&$A\sb{2}+A\sb{11}+D\sb{5}$&90&$2\;A\sb{1}+2\;A\sb{2}+A\sb{12}$&$23$\\ 
\hline 37&$A\sb{4}+2\;D\sb{7}$&213&$A\sb{4}+D\sb{14}$&$21$\\ 
\hline 38&$A\sb{3}+2\;A\sb{4}+D\sb{7}$&44&$2\;A\sb{1}+2\;A\sb{4}+A\sb{8}$&$23$\\ 
\hline 39&$2\;A\sb{2}+A\sb{3}+A\sb{4}+D\sb{7}$&20 & $A\sb{1}+2\;A\sb{2}+A\sb{3}+A\sb{4}+A\sb{6}$&$24$\\ 
\hline 40&$A\sb{2}+A\sb{4}+A\sb{5}+D\sb{7}$&23&$A\sb{1}+A\sb{2}+A\sb{4}+A\sb{5}+A\sb{6}$&$23$\\ 
\hline 41&$A\sb{1}+2\;A\sb{2}+A\sb{6}+D\sb{7}$&14&$2\;A\sb{1}+2\;A\sb{2}+2\;A\sb{6}$&$24$\\ 
\hline 
\end{SubConfigTable}\vfill\eject 
\begin{SubConfigTable} 
\hline 42&$2\;A\sb{2}+A\sb{7}+D\sb{7}$&90&$2\;A\sb{1}+2\;A\sb{2}+A\sb{12}$&$23$\\ 
\hline 43&$A\sb{4}+A\sb{7}+D\sb{7}$&44&$2\;A\sb{1}+2\;A\sb{4}+A\sb{8}$&$23 $\\ 
\hline 44&$A\sb{1}+A\sb{2}+A\sb{8}+D\sb{7}$&50&$2\;A\sb{1}+A\sb{2}+A\sb{6}+A\sb{8}$&$23$\\ 
\hline 45&$A\sb{3}+A\sb{8}+D\sb{7}$&44&$2\;A\sb{1}+2\;A\sb{4}+A\sb{8}$&$23$\\ 
\hline 46&$A\sb{11}+D\sb{7}$ & 320&$D\sb{10}+E\sb{8}$&$22$\\ 
\hline 47&$A\sb{2}+A\sb{4}+D\sb{5}+D\sb{7}$ & 200&$A\sb{2}+A\sb{5}+D\sb{11}$&$22$\\ 
\hline 48&$A\sb{6}+D\sb{5}+D\sb{7}$&186&$A\sb{9}+D\sb{9}$&$21$\\ 
\hline 49&$A\sb{2}+2\;A\sb{4}+D\sb{8}$&66&$A\sb{2}+A\sb{7}+A\sb{9}$&$ 21$\\ 
\hline 50&$A\sb{4}+A\sb{6}+D\sb{8}$&23&$A\sb{1}+A\sb{2}+A\sb{4}+A\sb{5}+A\sb{6}$&$23$\\ 
\hline 51&$A\sb{2}+A\sb{8}+D\sb{8}$&50&$2\;A\sb{1}+A\sb{2}+A\sb{6}+A\sb{8}$&$23$\\ 
\hline 52&$A\sb{10}+D\sb{8}$&320&$D\sb{10}+E\sb{8}$&$22$\\ 
\hline 53&$A\sb{1}+2\;A\sb{4}+D\sb{9}$&44&$2\;A\sb{1}+2\;A\sb{4}+A\sb{8}$&$23$\\ 
\hline 54&$A\sb{2}+A\sb{3}+A\sb{4}+D\sb{9}$&46&$A\sb{1}+A\sb{2}+A\sb{3}+A\sb{4}+A\sb{8}$&$23$\\ 
\hline 55&$A\sb{3}+A\sb{6}+D\sb{9}$&76&$2\;A\sb{1}+A\sb{6}+A\sb{10}$&$22$\\ 
\hline 56&$A\sb{2}+A\sb{7}+D\sb{9}$&50&$2\;A\sb{1}+A\sb{2}+A\sb{6}+A\sb{8}$&$23$\\ 
\hline 57&$2\;A\sb{2}+D\sb{5}+D\sb{9}$&210&$2\;A\sb{2}+D\sb{14}$&$22$\\ 
\hline 58&$A\sb{2}+D\sb{7}+D\sb{9}$&186&$A\sb{9}+D\sb{9}$&$21$\\ 
\hline 59&$2\;A\sb{2}+A\sb{4}+D\sb{10}$&72&$2\;A\sb{1}+A\sb{2}+A\sb{4}+A\sb{10}$&$ 23$\\ 
\hline 60&$A\sb{3}+A\sb{4}+D\sb{11}$&44&$2\;A\sb{1}+2\;A\sb{4}+A\sb{8}$&$23$\\ 
\hline 61&$A\sb{7}+D\sb{11}$&320&$D\sb{10}+E\sb{8}$&$22$\\ 
\hline 62&$A\sb{2}+D\sb{5}+D\sb{11}$&186&$A\sb{9}+D\sb{9}$&$21$\\ 
\hline 63&$D\sb{7}+D\sb{11}$ & 218&$D\sb{18}$&$20$\\ 
\hline 64&$A\sb{2}+A\sb{4}+D\sb{12}$&72&$2\;A\sb{1}+A\sb{2}+A\sb{4}+A\sb{10}$&$23$\\ 
\hline 65&$A\sb{6}+D\sb{12}$ & 320&$D\sb{10}+E\sb{8}$&$22$\\ 
\hline 66&$A\sb{1}+2\;A\sb{2}+D\sb{13}$&90&$2\;A\sb{1}+2\;A\sb{2}+A\sb{12}$&$23$\\ 
\hline 67&$A\sb{2}+A\sb{3}+D\sb{13}$&72&$2\;A\sb{1}+A\sb{2}+A\sb{4}+A\sb{10}$&$23$\\ 
\hline 68&$A\sb{3}+D\sb{15}$&320&$D\sb{10}+E\sb{8}$&$22$\\ 
\hline 69&$A\sb{2}+D\sb{16}$ & 320&$D\sb{10}+E\sb{8}$&$22$\\ 
\hline 70&$2\;A\sb{1}+A\sb{4}+2\;E\sb{6}$&303&$A\sb{1}+A\sb{4}+A\sb{5}+E\sb{8}$&$23$\\ 
\hline 71&$2\;A\sb{1}+A\sb{2}+2\;A\sb{4}+E\sb{6}$&23&$A\sb{1}+A\sb{2}+A\sb{4}+A\sb{5}+A\sb{6}$&$23$\\ 
\hline 72&$A\sb{2}+2\;A\sb{3}+A\sb{4}+E\sb{6}$&46&$A\sb{1}+A\sb{2}+A\sb{3}+A\sb{4}+A\sb{8}$&$23$\\ 
\hline 73&$2\;A\sb{6}+E\sb{6}$&37 &$ A\sb{1}+2\;A\sb{2}+A\sb{6}+A\sb{7}$&$23$\\ 
\hline 74&$2\;A\sb{3}+A\sb{6}+E\sb{6}$&41&$A\sb{5}+A\sb{6}+A\sb{7}$&$21$\\ 
\hline 75&$A\sb{2}+A\sb{3}+A\sb{7}+E\sb{6}$&37&$A\sb{1}+2\;A\sb{2}+A\sb{6}+A\sb{7}$&$23$\\ 
\hline 76&$2\;A\sb{4}+D\sb{4}+E\sb{6}$&182&$A\sb{4}+A\sb{5}+D\sb{9}$&$22$\\ 
\hline 77&$A\sb{2}+A\sb{6}+D\sb{4}+E\sb{6}$&183&$A\sb{1}+A\sb{2}+A\sb{6}+D\sb{9}$&$23$\\ 
\hline 78&$A\sb{8}+D\sb{4}+E\sb{6}$&186&$A\sb{9}+D\sb{9}$&$21$\\ 
\hline 79&$A\sb{1}+D\sb{5}+2\;E\sb{6}$&320&$D\sb{10}+E\sb{8}$&$22$\\ 
\hline 80&$A\sb{2}+2\;D\sb{5}+E\sb{6}$&320&$D\sb{10}+E\sb{8}$&$22$\\ 
\hline 81&$A\sb{1}+A\sb{2}+A\sb{4}+D\sb{5}+E\sb{6}$&193&$A\sb{2}+A\sb{6}+D\sb{10}$&$22$\\ 
\hline 82&$A\sb{2}+A\sb{3}+D\sb{7}+E\sb{6}$&200&$A\sb{2}+A\sb{5}+D\sb{11}$&$22$\\ 
\hline 
\end{SubConfigTable}\vfill\eject
\begin{SubConfigTable} 
\hline 83&$A\sb{5}+D\sb{7}+E\sb{6}$&320&$D\sb{10}+E\sb{8}$&$22$\\ 
\hline 84&$A\sb{2}+D\sb{10}+E\sb{6}$&193&$A\sb{2}+A\sb{6}+D\sb{10}$&$22$\\ 
\hline 85&$A\sb{1}+A\sb{2}+2\;A\sb{4}+E\sb{7}$& 17&$2\;A\sb{1}+A\sb{2}+2\;A\sb{4}+A\sb{6}$&$24$\\ 
\hline 86&$A\sb{3}+2\;A\sb{4}+E\sb{7}$&18&$A\sb{1}+A\sb{3}+2\;A\sb{4}+A\sb{6}$&$ 23$\\ 
\hline 87&$2\;A\sb{2}+D\sb{7}+E\sb{7}$&210&$2\;A\sb{2}+D\sb{14}$&$22$\\ 
\hline 88&$A\sb{2}+2\;A\sb{4}+E\sb{8}$&36&$A\sb{2}+A\sb{4}+A\sb{5}+A\sb{7}$&$22$\\ 
\hline 89&$2\;A\sb{1}+2\;A\sb{2}+A\sb{4}+E\sb{8}$&30&$2\;A\sb{2}+A\sb{3}+A\sb{4}+A\sb{7}$&$23$\\ 
\hline 90 &$ 2\;A\sb{3}+A\sb{4}+E\sb{8}$&24&$A\sb{3}+A\sb{4}+A\sb{5}+A\sb{6}$&$22$\\ 
\hline 91&$A\sb{3}+A\sb{7}+E\sb{8}$&46&$A\sb{1}+A\sb{2}+A\sb{3}+A\sb{4}+A\sb{8}$&$23$\\ 
\hline 92&$A\sb{2}+A\sb{4}+D\sb{4}+E\sb{8}$&182&$A\sb{4}+A\sb{5}+D\sb{9}$&$22$\\ 
\hline 93&$A\sb{6}+D\sb{4}+E\sb{8}$&186&$A\sb{9}+D\sb{9}$&$21$\\
\hline 94&$A\sb{1}+2\;A\sb{2}+D\sb{5}+E\sb{8}$&210&$2\;A\sb{2}+D\sb{14}$&$22$\\ 
\hline 95&$A\sb{2}+A\sb{3}+D\sb{5}+E\sb{8}$&198&$2\;A\sb{2}+A\sb{3}+D\sb{11}$&$23$ \\ 
\hline 96&$A\sb{3}+D\sb{7}+E\sb{8}$& 213&$A\sb{4}+D\sb{14}$&$21$\\ 
\hline 97&$A\sb{2}+D\sb{8}+E\sb{8}$& 210&$2\;A\sb{2}+D\sb{14}$&$22$\\ 
\hline 98&$2\;A\sb{1}+A\sb{2}+E\sb{6}+E\sb{8}$&320&$D\sb{10}+E\sb{8}$&$22$\\ 
\hline \end{SubConfigTable}\vfill\eject
\begin{EllTab} 
\hline  \hline  1  & $6\,A\sb{3}$ & $\Z/ (4)  \times  \Z/ (4)$ & $4$ & $0$ & $4$ \\
\hline  \hline  2  & $2\,A\sb{1}+4\,A\sb{4}$ & $\Z/ (5)$ & $10$ & $0$ & $10$ \\  
\hline  \hline  3  & $2\,A\sb{2}+2\,A\sb{3}+2\,A\sb{4}$ & $(0)$ & $60$ & $0$ & $60$ \\  
\hline  \hline  4  & $3\,A\sb{1}+3\,A\sb{5}$ & $\Z/ (2)  \times  \Z/ (6)$ & $2$ & $0$ & $6$ \\  
\hline  \hline  5  & $4\,A\sb{2}+2\,A\sb{5}$ & $\Z/ (3)  \times  \Z/ (3)$ & $6$ & $0$ & $6$ \\  
\hline  \hline  6  & $A\sb{3}+3\,A\sb{5}$ & $\Z/ (6)$ & $4$ & $0$ & $6$ \\  
\hline  \hline  7  & $2\,A\sb{1}+2\,A\sb{3}+2\,A\sb{5}$ & $\Z/ (2)  \times  \Z/ (2)$ & $12$ & $0$ & $12$ \\  
\hline  \hline  8  & $A\sb{1}+2\,A\sb{2}+A\sb{3}+2\,A\sb{5}$ & $\Z/ (6)$ & $6$ & $0$ & $12$ \\  
\hline  \hline  9  & $2\,A\sb{4}+2\,A\sb{5}$ & $(0)$ & $30$ & $0$ & $30$ \\  
\hline  \hline  10  & $2\,A\sb{2}+A\sb{4}+2\,A\sb{5}$ & $\Z/ (3)$ & $6$ & $0$ & $30$ \\  
\hline  \hline  11  & $A\sb{1}+A\sb{3}+A\sb{4}+2\,A\sb{5}$ & $\Z/ (2)$ & $12$ & $0$ & $30$ \\  
\hline  \hline  12  & $A\sb{1}+A\sb{2}+2\,A\sb{3}+A\sb{4}+A\sb{5}$ & $\Z/ (2)$ & $24$ & $12$ & $36$ \\  
\hline  \hline  13  & $3\,A\sb{6}$ & $\Z/ (7)$ & $2$ & $1$ & $4$ \\  
\hline  \hline  14  & $2\,A\sb{1}+2\,A\sb{2}+2\,A\sb{6}$ & $(0)$ & $42$ & $0$ & $42$ \\  
\hline  \hline  15  & $2\,A\sb{3}+2\,A\sb{6}$ & $(0)$ & $28$ & $0$ & $28$ \\  
\hline  \hline  16  & $A\sb{2 }+A\sb{4}+2\,A\sb{6}$ & $(0)$ & $28$ & $7$ & $28$ \\ 
\hline  \hline  17  & $2\,A\sb{1}+A\sb{2}+2\,A\sb{4}+A\sb{6}$ & $(0)$ & $50$ & $20$ & $50$ \\  
\hline  \hline  18  & $A\sb{1}+A\sb{3}+2\,A\sb{4}+A\sb{6}$ & $(0)$ & $10$ & $0$ & $140$ \\  
\cline{4-6}  &  &  & $20$ & $0$ & $70$ \\  
\hline  \hline  19  & $A\sb{2}+2\,A\sb{3}+A\sb{4}+A\sb{6}$ & $(0)$ & $24$ & $12$ & $76$ \\  
\hline  \hline  20  & $A\sb{1}+2\,A\sb{2}+A\sb{3}+A\sb{4}+A\sb{6}$ & $(0)$ & $30$ & $0$ & $84$ \\  
\hline  \hline  21  & $2\,A\sb{1}+2\,A\sb{5}+A\sb{6}$ & $\Z/ (2)$ & $12$ & $6$ & $24$ \\  
\hline  \hline  22  & $A\sb{1}+2\,A\sb{3}+A\sb{5}+A\sb{6}$ & $\Z/ (2)$ & $4$ & $0$ & $84$ \\  
\hline  \hline  23  & $A\sb{1}+A\sb{2}+A\sb{4}+A\sb{5}+A\sb{6}$ & $(0)$ & $30$ & $0$ & $42$ \\  
\cline{4-6}  &  &  & $18$ & $6$ & $72$ \\  
\hline  \hline  24  & $A\sb{3}+A\sb{4}+A\sb{5}+A\sb{6}$ & $(0)$ & $12$ & $0$ & $70$ \\  
\hline  \hline  25  & $4\,A\sb{1}+2\,A\sb{7}$ & $\Z/ (2)  \times  \Z/ (4)$ & $4$ & $0$ & $4$ \\  
\hline  \hline  26  & $2\,A\sb{2}+2\,A\sb{7}$ & $(0)$ & $24$ & $0$ & $24$ \\  
\cline{3-6}  &  & $\Z/ (2)$ & $12$ & $0$ & $12$ \\  
\hline  \hline  27  & $A\sb{1}+A\sb{3}+2\,A\sb{7}$ & $\Z/ (8)$ & $2$ & $0$ & $4$ \\  
\hline  \hline  28  & $2\,A\sb{1}+3\,A\sb{3}+A\sb{7}$ & $\Z/ (2)  \times  \Z/ (4)$ & $4$ & $0$ & $8$ \\  
\hline  \hline  29  & $A\sb{2}+3\,A\sb{3}+A\sb{7}$ & $\Z/ (4)$ & $4$ & $0$ & $24$ \\  
\hline  \hline  30  & $2\,A\sb{2}+A\sb{3}+A\sb{4}+A\sb{7}$ & $(0)$ & $12$ & $0$ & $120$ \\  
\hline  \hline  31  & $2\,A\sb{1}+A\sb{2}+A\sb{3}+A\sb{4}+A\sb{7}$ & $\Z/ (2)$ & $20$ & $0$ & $24$ \\  
\hline  \hline  32  & $A\sb{1}+2\,A\sb{5}+A\sb{7}$ & $\Z/ (2)$ & $6$ & $0$ & $24$ \\  
\hline  \hline  33  & $3\,A\sb{1}+A\sb{3}+A\sb{5}+A\sb{7}$ & $\Z/ (2)  \times  \Z/ (2)$ & $8$ & $0$ & $12$ \\  
\hline  \end{EllTab}  \vfill  \eject  
\begin{EllTab} 
\hline  \hline  34  & $A\sb{1}+A\sb{2}+A\sb{3}+A\sb{5}+A\sb{7}$ & $\Z/ (2)$ & $12$ & $0$ & $24$
\\  
\hline  \hline  35  & $2\,A\sb{1}+A\sb{4}+A\sb{5}+A\sb{7}$ & $\Z/ (2)$ & $2$ & $0$ & $120$ \\  
\hline  \hline  36  & $A\sb{2}+A\sb{4}+A\sb{5}+A\sb{7}$ & $(0)$ & $6$ & $0$ & $120$ \\  
\cline{4-6}  &  &  & $24$ & $0$ & $30$ \\  
\hline  \hline  37  & $A\sb{1}+2\,A\sb{2}+A\sb{6}+A\sb{7}$ & $(0)$ & $24$ & $0$ & $42$ \\  
\hline  \hline  38  & $2\,A\sb{1}+A\sb{3}+A\sb{6}+A\sb{7}$ & $\Z/ (2)$ & $12$ & $4$ & $20$ \\  
\hline  \hline  39  & $A\sb{2}+A\sb{3}+A\sb{6}+A\sb{7}$ & $(0)$ & $4$ & $0$ & $168$ \\  
\hline  \hline  40  & $A\sb{1}+A\sb{4}+A\sb{6}+A\sb{7}$ & $(0)$ & $2$ & $0$ & $280$ \\  
\cline{4-6}  &  &  & $18$ & $4$ & $32$ \\  
\hline  \hline  41  & $A\sb{5}+A\sb{6}+A\sb{7}$ & $(0)$ & $16$ & $4$ & $22$ \\  
\hline  \hline  42  & $2\,A\sb{1}+2\,A\sb{8}$ & $(0)$ & $18$ & $0$ & $18$ \\  
\cline{3-6}  &  & $\Z/ (3)$ & $4$ & $2$ & $10$ \\  
\hline  \hline  43  & $A\sb{1}+3\,A\sb{2}+A\sb{3}+A\sb{8}$ & $\Z/ (3)$ & $12$ & $0$ & $18$ \\  
\hline  \hline  44  & $2\,A\sb{1}+2\,A\sb{4}+A\sb{8}$ & $(0)$ & $20$ & $10$ & $50$ \\  
\hline  \hline  45  & $3\,A\sb{2}+A\sb{4}+A\sb{8}$ & $\Z/ (3)$ & $12$ & $3$ & $12$ \\  
\hline  \hline  46  & $A\sb{1}+A\sb{2}+A\sb{3}+A\sb{4}+A\sb{8}$ & $(0)$ & $6$ & $0$ & $180$ \\  
\hline  \hline  47  & $A\sb{1}+2\,A\sb{2}+A\sb{5}+A\sb{8}$ & $\Z/ (3)$ & $6$ & $0$ & $18$ \\  
\hline  \hline  48  & $A\sb{2}+A\sb{3}+A\sb{5}+A\sb{8}$ & $\Z/ (3)$ & $4$ & $0$ & $18$ \\  
\hline  \hline  49  & $A\sb{1}+A\sb{4}+A\sb{5}+A\sb{8}$ & $(0)$ & $18$ & $0$ & $30$ \\  
\hline  \hline  50  & $2\,A\sb{1}+A\sb{2}+A\sb{6}+A\sb{8}$ & $(0)$ & $18$ & $0$ & $42$ \\  
\hline  \hline  51  & $A\sb{1}+A\sb{3}+A\sb{6}+A\sb{8}$ & $(0)$ & $10$ & $4$ & $52$ \\  
\hline  \hline  52  & $A\sb{4}+A\sb{6}+A\sb{8}$ & $(0)$ & $18$ & $9$ & $22$ \\  
\hline  \hline  53  & $A\sb{1}+A\sb{2}+A\sb{7}+A\sb{8}$ & $(0)$ & $18$ & $0$ & $24$ \\  
\hline  \hline  54  & $2\,A\sb{9}$ & $(0)$ & $10$ & $0$ & $10$ \\  
\cline{3-6}  &  & $\Z/ (5)$ & $2$ & $0$ & $2$ \\  
\hline  \hline  55  & $A\sb{1}+A\sb{2}+2\,A\sb{3}+A\sb{9}$ & $\Z/ (2)$ & $4$ & $0$ & $60$ \\  
\hline  \hline  56  & $2\,A\sb{1}+2\,A\sb{2}+A\sb{3}+A\sb{9}$ & $\Z/ (2)$ & $6$ & $0$ & $60$ \\  
\hline  \hline  57  & $A\sb{1}+2\,A\sb{4}+A\sb{9}$ & $\Z/ (5)$ & $2$ & $0$ & $10$ \\  
\hline  \hline  58  & $3\,A\sb{1}+A\sb{2}+A\sb{4}+A\sb{9}$ & $\Z/ (2)$ & $20$ & $10$ & $20$ \\  
\hline  \hline  59  & $2\,A\sb{1}+A\sb{3}+A\sb{4}+A\sb{9}$ & $\Z/ (2)$ & $10$ & $0$ & $20$ \\  
\hline  \hline  60  & $2\,A\sb{1}+A\sb{2}+A\sb{5}+A\sb{9}$ & $\Z/ (2)$ & $12$ & $6$ & $18$ \\  
\hline  \hline  61  & $A\sb{1}+A\sb{3}+A\sb{5}+A\sb{9}$ & $\Z/ (2)$ & $10$ & $0$ & $12$ \\  
\hline  \hline  62  & $A\sb{4}+A\sb{5}+A\sb{9}$ & $(0)$ & $10$ & $0$ & $30$ \\  
\cline{3-6}  &  & $\Z/ (2)$ & $10$ & $5$ & $10$ \\ 
\hline  \hline  63   & $3\,A\sb{1}+A\sb{6}+A\sb{9}$ & $\Z/ (2)$ & $4$ & $2$ & $36$ \\  
\hline  \hline  64  & $A\sb{1}+A\sb{2}+A\sb{6}+A\sb{9}$ & $(0)$ & $10$ & $0$ & $42$ \\  
\hline  \end{EllTab}  \vfill  \eject  
\begin{EllTab} 
\hline  \hline  65  & $A\sb{3}+A\sb{6}+A\sb{9}$ & $(0)$ & $2$ & $0$ & $140$ \\  
\hline  \hline  66  & $A\sb{2}+A\sb{7}+A\sb{9}$ & $(0)$ & $10$ & $0$ & $24$ \\  
\hline  \hline  67  & $A\sb{1}+A\sb{8}+A\sb{9}$ & $(0)$ & $10$ & $0$ & $18$ \\  
\hline  \hline  68  & $A\sb{2}+2\,A\sb{3}+A\sb{10}$ & $(0)$ & $24$ & $12$ & $28$ \\  
\hline  \hline  69  & $A\sb{1}+2\,A\sb{2}+A\sb{3}+A\sb{10}$ & $(0)$ & $12$ & $0$ & $66$ \\  
\hline  \hline  70  & $2\,A\sb{4}+A\sb{10}$ & $(0)$ & $10$ & $5$ & $30$ \\  
\hline  \hline  71  & $2\,A\sb{2}+A\sb{4}+A\sb{10}$ & $(0)$ & $6$ & $3$ & $84$ \\  
\cline{4-6} &  &  & $24$ & $9$ & $24$ \\  
\hline  \hline  72  & $2\,A\sb{1}+A\sb{2}+A\sb{4}+A\sb{10}$ & $(0)$ & $2$ & $0$ & $330$ \\  
\hline  \hline  73  & $A\sb{1}+A\sb{3}+A\sb{4}+A\sb{10}$ & $(0)$ & $20$ & $0$ & $22$ \\  
\cline{4-6}   &  &  & $12$ & $4$ & $38$ \\  
\hline  \hline  74  & $A\sb{1}+A\sb{2}+A\sb{5}+A\sb{10}$ & $(0)$ & $6$ & $0$ & $66$ \\  
\cline{4-6}  & &  & $18$ & $6$ & $24$ \\  
\hline  \hline  75  & $A\sb{3}+A\sb{5}+A\sb{10}$ & $(0)$ & $4$ & $0$ & $66$ \\  
\cline{4-6}  &  &  & $12$ & $0$ & $22$ \\  
\hline  \hline  76  & $2\,A\sb{1}+A\sb{6}+A\sb{10}$ & $(0)$ & $12$ & $2$ & $26$ \\  
\hline  \hline  77  & $A\sb{2}+A\sb{6}+A\sb{10}$ & $(0)$ & $4$ & $1$ & $58$ \\  
\cline{4-6}  &  &  & $16$ & $5$ & $16$ \\  
\hline  \hline  78  & $A\sb{1 }+A\sb{7}+A\sb{10}$ & $(0)$ & $2$ & $0$ & $88$ \\  
\cline{4-6}  &  &  & $10$ & $2$ & $18$ \\  
\hline  \hline  79  & $A\sb{8}+A\sb{10}$ & $(0)$ & $10$ & $1$ & $10$ \\  
\hline  \hline  80  & $A\sb{1}+3\,A\sb{2}+A\sb{11}$ & $\Z/ (3)$ & $6$ & $0$ & $12$ \\  
\hline  \hline  81  & $3\,A\sb{1}+2\,A\sb{2}+A\sb{11}$ & $\Z/ (6)$ & $2$ & $0$ & $12$ \\  
\hline  \hline  82  & $A\sb{1}+2\,A\sb{3}+A\sb{11}$ & $\Z/ (4)$ & $4$ & $0$ & $6$ \\  
\hline  \hline  83  & $2\,A\sb{2} +A\sb{3}+A\sb{11}$ & $\Z/ (3)$ & $4$ & $0$ & $12$ \\  
\cline{3-6}  &  & $\Z/ (6)$ & $4$ & $2$ & $4$  \\  
\hline  \hline  84  & $2\,A\sb{1}+A\sb{2}+A\sb{3}+A\sb{11}$ & $\Z/ (4)$ & $6$ & $0$ & $6$ \\  
\cline{3-6}  &  & $\Z/ (2)$ & $12$ & $0$ & $12$ \\  
\hline  \hline  85  & $3\,A\sb{1}+A\sb{4}+A\sb{11}$ & $\Z/ (2)$ & $6$ & $0$ & $20$ \\  
\hline  \hline  86  & $A\sb{1}+A\sb{2}+A\sb{4}+A\sb{11}$ & $(0)$ & $12$ & $0$ & $30$ \\  
\hline  \hline  87  & $2\,A\sb{1}+A\sb{5}+A\sb{11}$ & $\Z/ (2)$ & $6$ & $0$ & $12$ \\  
\cline{3-6}  &  & $\Z/ (6)$ & $2$ & $0$ & $4$ \\  
\hline  \hline  88  & $A\sb{2}+A\sb{5}+A\sb{11}$ & $\Z/ (3)$ & $4$ & $0$ & $6$ \\  
\hline  \hline  89  & $A\sb{1}+A\sb{6}+A\sb{11}$ & $(0)$ & $4$ & $0$ & $42$ \\  
\hline  \hline  90  & $2\,A\sb{1}+2\,A\sb{2}+A\sb{12}$ & $(0)$ & $12$ & $6$ & $42$ \\  
\hline  \hline  91  & $A\sb{1}+A\sb{2}+A\sb{3}+A\sb{12}$ & $(0)$ & $6$ & $0$ & $52$ \\  \hline  
\end{EllTab} 
\vfill  \eject  
\begin{EllTab} 
\hline  \hline  92  & $2\,A\sb{1}+A\sb{4}+A\sb{12}$ & $(0)$ & $2$ & $0$ & $130$ \\  
\cline{4-6}  &  &  & $18$ & $8$ & $18$ \\  
\hline  \hline  93  & $A\sb{2}+A\sb{4}+A\sb{12}$ & $(0)$ & $6$ & $3$ & $34$ \\  
\hline  \hline  94  & $A\sb{1}+A\sb{5}+A\sb{12}$ & $(0)$ & $10$ & $2$ & $16$ \\  
\hline  \hline  95  & $A\sb{6}+A\sb{12}$ & $(0)$ & $2$ & $1$ & $46$ \\  
\hline  \hline  96  & $A\sb{1}+2\,A\sb{2}+A\sb{13}$ & $(0)$ & $6$ & $0$ & $42$ \\  
\cline{3-6}  &  & $\Z/ (2)$ & $6$ & $3$ & $12$ \\  
\hline  \hline  97  & $3\,A\sb{1}+A\sb{2}+A\sb{13}$ & $\Z/ (2)$ & $2$ & $0$ & $42$ \\  
\hline  \hline  98  & $2\,A\sb{1}+A\sb{3}+A\sb{13}$ & $\Z/ (2)$ & $6$ & $2$ & $10$ \\  
\hline  \hline  99  & $A\sb{2}+A\sb{3}+A\sb{13}$ & $(0)$ & $4$ & $0$ & $42$ \\  
\hline  \hline  100  & $A\sb{1}+A\sb{4}+A\sb{13}$ & $(0)$ & $2$ & $0$ & $70$ \\  
\cline{4-6} &  &  & $8$ & $2$ & $18$ \\  
\cline{3-6} &  & $\Z/ (2)$ & $2$ & $1$ & $18$ \\  
\hline  \hline  101  & $A\sb{5}+A\sb{13}$ & $(0)$ & $4$ & $2$ & $22$ \\  
\hline  \hline  102  & $2\,A\sb{2}+A\sb{14}$ & $\Z/ (3)$ & $4$ & $1$ & $4$ \\  
\hline  \hline  103  & $2\,A\sb{1}+A\sb{2}+A\sb{14}$ & $(0)$ & $12$ & $6$ & $18$ \\  
\cline{3-6}  &  & $\Z/ (3)$ & $2$ & $0$ & $10$ \\ 
\hline  \hline  104  & $A\sb{1}+A\sb{3}+A\sb{14}$ & $(0)$ & $10$ & $0$ & $12$ \\  
\hline  \hline  105  & $A\sb{4}+A\sb{14}$ & $(0)$ & $10$ & $5$ & $10$ \\  
\hline  \hline  106  & $3\, A\sb{1}+A\sb{15}$ & $\Z/ (4)$ & $2$ & $0$ & $4$ \\  
\hline  \hline  107  & $A\sb{1}+ A\sb{2}+A\sb{15}$ & $(0)$ & $10$ & $2$ & $10$ \\  
\cline{3-6}  &  & $\Z/ (2)$ & $4$ & $0$ & $6$ \\  
\hline  \hline  108  & $A\sb{3}+A\sb{15}$ & $\Z/ (4)$ & $2$ & $0$ & $2$ \\  
\hline  \hline  109  & $2\,A\sb{1}+A\sb{16}$ & $(0)$ & $2$ & $0$ & $34$ \\  
\cline{4-6}   &  &  & $4$ & $2$ & $18$ \\  
\hline  \hline  110  & $A\sb{2}+A\sb{16}$ & $(0)$ & $6$ & $3$ & $10$ \\  
\hline  \hline  111  & $A\sb{1}+A\sb{17}$ & $(0)$ & $4$ & $2$ & $10$ \\  
\cline{3-6}  &  & $\Z/ (3)$ & $2$ & $0$ & $2$ \\  
\hline  \hline  112  & $A\sb{18}$ & $(0)$ & $2$ & $1$ & $10$ \\  
\hline  \hline  113  & $2\,A\sb{4}+2\,D\sb{5}$ & $(0)$ & $20$ & $0$ & $20$ \\  
\hline  \hline  114  & $A\sb{3}+2\,A\sb{5}+D\sb{5}$ & $\Z/ (2)$ & $12$ & $0$ & $12$ \\  
\hline  \hline  115  & $2\,A\sb{4}+A\sb{5}+D\sb{5}$ & $(0)$ & $20$ & $0$ & $30$ \\  
\hline  \hline  116  & $A\sb{1}+A\sb{3}+A\sb{4}+A\sb{5}+D\sb{5}$ & $\Z/ (2)$ & $12$ & $0$ & $20$ \\  
\hline  \hline  117  & $A\sb{1}+2\,A\sb{6}+D\sb{5}$ & $(0)$ & $14$ & $0$ & $28$ \\  
\hline  \hline  118  & $2\,A\sb{2}+A\sb{3}+A\sb{6}+D\sb{5}$ & $(0)$ & $12$ & $0$ & $84$ \\  
\hline  \hline  119  & $A\sb{1}+A\sb{2}+A\sb{4}+A\sb{6}+D\sb{5}$ & $(0)$ & $20$ & $0$ & $42$ \\  
\hline  \end{EllTab}  \vfill  \eject  
\begin{EllTab} 
\hline  \hline  120  & $A\sb{2}+A\sb{5}+A\sb{6}+D\sb{5}$ & $(0)$ & $6$ & $0$ & $84$ \\  
\cline{4-6}  &  &  & $12$ & $0$ & $42$ \\
\hline  \hline  121  & $A\sb{1}+A\sb{7}+2\,D\sb{5}$ & $\Z/ (4)$ & $2$ & $0$ & $8$ \\
\hline  \hline  122  & $A\sb{1}+A\sb{2}+A\sb{3}+A\sb{7}+D\sb{5}$ & $\Z/ (4)$ & $6$ & $0$ & $8$ \\  
\hline  \hline  123  & $2\,A\sb{1}+A\sb{4}+A\sb{7}+D\sb{5}$ & $\Z/ (2)$ & $8$ & $0$ & $20$ \\
\hline  \hline  124  & $A\sb{8}+2\,D\sb{5}$ & $(0)$ & $8$ & $4$ & $20$ \\  
\hline  \hline  125  & $A\sb{1}+A\sb{4}+A\sb{8}+D\sb{5}$ & $(0)$ & $2$ & $0$ & $180$ \\
\cline{4-6}  &  &  & $18$ & $0$ & $20$ \\
\hline  \hline  126  & $A\sb{5}+A\sb{8}+D\sb{5}$ & $(0)$ & $12$ & $0$ & $18$ \\
\hline  \hline  127  & $2\,A\sb{2}+A\sb{9}+D\sb{5}$ & $(0)$ & $6$ & $0$ & $60$ \\  
\hline  \hline  128  & $2\,A\sb{1}+A\sb{2}+ A\sb{9}+D\sb{5}$ & $\Z/ (2)$ & $2$ & $0$ & $60$ \\  
\hline  \hline  129  & $A\sb{1}+A\sb{3}+A\sb{9}+D\sb{5}$ & $\Z/ (2)$ & $8$ & $4$ & $12$ \\  
\hline  \hline  130  & $A\sb{4}+A\sb{9}+D\sb{5}$ & $(0)$ & $10$ & $0$ & $20$ \\  
\hline  \hline  131  & $A\sb{1}+A\sb{2}+A\sb{10}+D\sb{5}$ & $(0)$ & $14$ & $4$ & $20$ \\  
\hline  \hline  132  & $2\, A\sb{1}+A\sb{11}+D\sb{5}$& $\Z/ (4)$ & $2$ & $0$ & $6$ \\  
\hline  \hline  133  & $A\sb{2}+A\sb{11}+D\sb{5}$ & $\Z/  (2)$ & $6$ & $0$ & $6$ \\  
\hline  \hline  134   & $A\sb{1}+A\sb{12}+D\sb{5}$ & $(0)$ & $2$ & $0$ & $52$ \\  
\cline{4-6}  &  &  & $6$ & $2$ & $18$ \\  
\hline  \hline  135  & $A\sb{13}+D\sb{5}$ & $(0)$ & $6$ & $2$ & $10$ \\  
\hline  \hline  136  & $3\,D\sb{6}$ & $\Z/ (2)  \times  \Z/ (2)$ & $2$ & $0$ & $2$ \\  
\hline  \hline  137  & $2\,A\sb{3}+2\,D\sb{6}$ & $\Z/ (2)  \times  \Z/ (2)$ & $4$ & $0$ & $4$ \\  
\hline  \hline  138  & $2\,A\sb{2}+2\,A\sb{4}+D\sb{6}$ & $(0)$ & $30$ & $0$ & $30$ \\  
\hline  \hline  139  & $2\,A\sb{1}+2\,A\sb{5}+D\sb{6}$ & $\Z/ (2)  \times  \Z/ (2)$ & $6$ & $0$ & $6$ \\  
\hline  \hline  140  & $A\sb{1}+2\,A\sb{3}+A\sb{5}+D\sb{6}$ & $\Z/ (2)  \times  \Z/ (2)$ & $4$ & $0$ & $12$ \\  
\hline  \hline  141  & $A\sb{3}+A\sb{4}+A\sb{5}+D\sb{6}$ & $\Z/ (2)$ & $4$ & $0$ & $30$ \\  
\hline  \hline  142  & $2\,A\sb{6}+D\sb{6}$ & $(0)$ & $14$ & $0$ & $14$ \\  
\hline  \hline  143  & $A\sb{2}+A\sb{4}+A\sb{6}+D\sb{6}$ & $(0)$ & $6$ & $0$ & $70$ \\  
\hline  \hline  144  & $A\sb{1}+2\,A\sb{2}+A\sb{7}+D\sb{6}$ & $\Z/ (2)$ & $6$ & $0$ & $24$ \\  
\hline  \hline  145  & $A\sb{2}+A\sb{3}+A\sb{7}+D\sb{6}$ & $\Z/ (2)$ & $4$ & $0$ & $24$ \\  
\hline  \hline  146  & $A\sb{1}+A\sb{4}+A\sb{7}+D\sb{6}$ & $\Z/ (2)$ & $6$ & $2$ & $14$ \\  
\hline  \hline  147  & $A\sb{4}+A\sb{8}+D\sb{6}$ & $(0)$ & $4$ & $2$ & $46$ \\  
\hline  \hline  148  & $A\sb{1}+A\sb{2}+A\sb{9}+D\sb{6}$ & $\Z/ (2)$ & $6$ & $0$ & $10$ \\  
\cline{3-6}  &  & $\Z/ (2)$ & $4$ & $2$ & $16$ \\  
\hline  \hline  149  & $A\sb{3}+A\sb{9}+D\sb{6}$ & $\Z/ (2)$ & $4$ & $0$ & $10$ \\  
\hline  \hline  150  & $A\sb{2}+A\sb{10}+D\sb{6}$ & $(0)$ & $6$ & $0$ & $22$ \\  
\hline  \hline  151  & $A\sb{1}+A\sb{11}+D\sb{6}$ & $\Z/ (2)$ & $4$ & $0$ & $6$ \\  
\hline  \end{EllTab}  \vfill  \eject  
\begin{EllTab} 
\hline  \hline  152  & $A\sb{12}+D\sb{6}$ & $(0)$ & $4$ & $2$ & $14$ \\  
\hline  \hline  153  & $A\sb{2}+A\sb{5}+D\sb{5}+D\sb{6}$ & $\Z/ (2)$ & $6$ & $0$ & $12$ \\
\hline  \hline  154  & $A\sb{7}+D\sb{5}+D\sb{6}$ & $\Z/ (2)$ & $4$ & $0$ & $8$ \\  
\hline  \hline  155  & $2\,A\sb{2}+2\,D\sb{7}$ & $(0)$ & $12$ & $0$ & $12$ \\  
\hline  \hline  156  & $A\sb{2}+3\,A\sb{3}+D\sb{7}$ & $\Z/ (4)$ & $8$ & $4$ & $8$ \\  
\hline  \hline  157  & $A\sb{1}+A\sb{2}+2\,A\sb{4}+D\sb{7}$ & $(0)$ & $10$ & $0$ & $60$ \\  
\hline  \hline  158  & $A\sb{2}+A\sb{3}+A\sb{6}+D\sb{7}$ & $(0)$ & $8$ & $4$ & $44$ \\  
\hline  \hline  159  & $A\sb{1}+A\sb{4}+A\sb{6}+D\sb{7}$ & $(0)$ & $4$ & $0$ & $70$ \\  
\hline  \hline  160  & $A\sb{5}+A\sb{6}+D\sb{7}$ & $(0)$ & $2$ & $0$ & $84$ \\  
\hline  \hline  161  & $2\,A\sb{1}+A\sb{2}+A\sb{7}+D\sb{7}$ & $\Z/ (2)$ & $4$ & $0$ & $24$ \\  
\hline  \hline  162  & $A\sb{1}+A\sb{3}+A\sb{7}+D\sb{7}$ & $\Z/ (4)$ & $2$ & $0$ & $8$ \\  
\hline  \hline  163  & $2\,A\sb{1}+A\sb{9}+D\sb{7}$ & $\Z/ (2)$ & $4$ & $0$ & $10$ \\  
\hline  \hline  164  & $A\sb{2}+A\sb{9}+D\sb{7}$ & $(0)$ & $2$ & $0$ & $60$ \\  
\hline  \hline  165  & $A\sb{1}+A\sb{10}+D\sb{7}$ & $(0)$ & $4$ & $0$ & $22$ \\  
\hline  \hline  166  & $A\sb{11}+D\sb{7}$ & $\Z/ (4)$ & $2$ & $1$ & $2$ \\  
\hline  \hline  167  & $A\sb{1}+A\sb{5}+D\sb{5}+D\sb{7}$ & $\Z/ (2)$ & $4$ & $0$ & $12$ \\  
\hline  \hline  168  & $A\sb{5}+D\sb{6}+D\sb{7}$ & $\Z/ (2)$ & $2$ & $0$ & $12$ \\  
\hline  \hline  169  & $2\,A\sb{1}+2\,D\sb{8}$ & $\Z/ (2)  \times  \Z/ (2)$ & $2$ & $0$ & $2$ \\  
\hline  \hline  170  & $2\,A\sb{2}+2\,A\sb{3}+D\sb{8}$ & $\Z/ (2)$ & $12$ & $0$ & $12$ \\  
\hline  \hline  171  & $2\,A\sb{5}+D\sb{8}$ & $\Z/ (2)$ & $6$ & $0$ & $6$ \\  
\hline  \hline  172  & $2\,A\sb{1}+A\sb{3}+A\sb{5}+D\sb{8}$ & $\Z/ (2)  \times  \Z/ (2)$ & $2$ & $0$ & $12$ \\  
\hline  \hline  173  & $A\sb{1}+A\sb{4}+A\sb{5}+D\sb{8}$ & $\Z/ (2)$ & $2$ & $0$ & $30$ \\  
\hline  \hline  174  & $2\,A\sb{2}+A\sb{6}+D\sb{8}$ & $(0)$ & $12$ & $6$ & $24$ \\  
\hline  \hline  175  & $A\sb{1}+A\sb{2}+A\sb{7}+D\sb{8}$ & $\Z/ (2)$ & $2$ & $0$ & $24$ \\  
\hline  \hline  176  & $A\sb{1}+A\sb{9}+D\sb{8}$ & $\Z/ (2)$ & $2$ & $0$ & $10$ \\  
\hline  \hline  177  & $2\,D\sb{5}+D\sb{8}$ & $\Z/ (2)$ & $4$ & $0$ & $4$ \\  
\hline  \hline  178  & $A\sb{1}+A\sb{3}+D\sb{6}+D\sb{8}$ & $\Z/ (2)  \times  \Z/ (2)$ & $2$ & $0$ & $4$ \\  
\hline  \hline  179  & $2\,D\sb{9}$ & $(0)$ & $4$ & $0$ & $4$ \\  
\hline  \hline  180  & $A\sb{1}+2\,A\sb{2}+A\sb{4}+D\sb{9}$ & $(0)$ & $12$ & $0$ & $30$ \\  
\hline  \hline  181  & $A\sb{1}+A\sb{3}+A\sb{5}+D\sb{9}$ & $\Z/ (2)$ & $4$ & $0$ & $12$ \\  
\hline  \hline  182  & $A\sb{4}+A\sb{5}+D\sb{9}$ & $(0)$ & $4$ & $0$ & $30$ \\  
\hline  \hline  183  & $A\sb{1}+A\sb{2}+A\sb{6}+D\sb{9}$ & $(0)$ & $4$ & $0$ & $42$ \\  
\hline  \hline  184  & $2\,A\sb{1}+A\sb{7}+D\sb{9}$ & $\Z/ (2)$ & $4$ & $0$ & $8$ \\
\hline  \hline  185  & $A\sb{1}+A\sb{8}+D\sb{9}$ & $(0)$ & $4$ & $0$ & $18$ \\  
\hline  \hline  186  & $A\sb{9}+D\sb{9}$ & $(0)$ & $4$ & $0$ & $10$ \\
\hline  \hline  187  & $A\sb{4}+D\sb{5}+D\sb{9}$ & $(0)$ & $4$ & $0$ & $20$ \\  
\hline  \end{EllTab}  \vfill  \eject  
\begin{EllTab} 
\hline  \hline  188  & $2\,A\sb{1}+2\,A\sb{3}+D\sb{10}$ & $\Z/ (2)  \times  \Z/ (2)$ & $4$ & $0$
& $4$ \\  
\hline  \hline  189  & $2\,A\sb{4}+D\sb{10}$ & $(0)$ & $10$ & $0$ & $10$ \\  
\hline  \hline  190  & $A\sb{1}+A\sb{3}+A\sb{4}+D\sb{10}$ & $\Z/ (2)$ & $2$ & $0$ & $20$ \\  
\hline  \hline  191  & $3\,A\sb{1}+A\sb{5}+D\sb{10}$ & $\Z/ (2)  \times  \Z/ (2)$ & $4$ & $2$ & $4$ \\  
\hline  \hline  192  & $A\sb{3}+A\sb{5}+D\sb{10}$ & $\Z/ (2)$ & $2$ & $0$ & $12$ \\  
\hline  \hline  193  & $A\sb{2}+A\sb{6}+D\sb{10}$ & $(0)$ & $2$ & $0$ & $42$ \\  
\hline  \hline  194  & $A\sb{8}+D\sb{10}$ & $(0)$ & $2$ & $0$ & $18$ \\  
\hline  \hline  195  & $A\sb{1}+A\sb{2}+D\sb{5}+D\sb{10}$ & $\Z/ (2)$ & $4$ & $0$ & $6$ \\  
\hline  \hline  196  & $A\sb{2}+D\sb{6}+D\sb{10}$ & $\Z/ (2)$ & $2$ & $0$ & $6$ \\  
\hline  \hline  197  & $A\sb{1}+D\sb{7}+D\sb{10}$ & $\Z/ (2)$ & $2$ & $0$ & $4$ \\  
\hline  \hline  198  & $2\,A\sb{2}+A\sb{3}+D\sb{ 11}$ & $(0)$ & $12$ & $0$ & $12$ \\  
\hline  \hline  199  & $A\sb{1}+A\sb{2}+A\sb{4}+D\sb{ 11}$ & $(0)$ & $6$ & $0$ & $20$ \\  
\hline  \hline  200 & $A\sb{2}+A\sb{5}+D\sb{11}$ & $(0)$ & $6$ & $0$ & $12$ \\  
\hline  \hline  201  & $A\sb{1}+A\sb{6}+D\sb{11}$ & $(0)$ & $6$ & $2$ & $10$ \\  
\hline  \hline  202  & $2\,A\sb{1}+2\,A\sb{2}+D\sb{12}$ & $\Z/ (2)$ & $6$ & $0$ & $6$ \\  
\hline  \hline  203  & $A\sb{1}+A\sb{2}+A\sb{3}+D\sb{12}$ & $\Z/ (2)$ & $4$ & $0$ & $6$ \\  
\hline  \hline  204  & $2\,A\sb{1}+A\sb{4}+D\sb{ 12}$ & $\Z/ (2)$ & $4$ & $2$ & $6$ \\  
\hline  \hline  205  & $A\sb{1}+D\sb{5}+D\sb{12}$ & $\Z/ (2)$ & $2$ & $0$ & $4$ \\  
\hline  \hline  206  & $D\sb{6}+D\sb{12}$ & $\Z/ (2)$ & $2$ & $0$ & $2$ \\  
\hline  \hline  207  & $A\sb{1}+A\sb{4}+D\sb{13}$ & $(0)$ & $2$ & $0$ & $20$ \\  
\hline  \hline  208  & $A\sb{5}+D\sb{13}$ & $(0)$ & $2$ & $0$ & $12$ \\  
\hline  \hline  209  & $D\sb{5}+D\sb{13}$& $(0)$ & $4$ & $0$ & $4$ \\  
\hline  \hline  210  & $2\,A\sb{2}+D\sb{14}$ & $(0)$ & $6$ & $0$ & $6$ \\  
\hline  \hline  211  & $2\,A\sb{1}+A\sb{2}+D\sb{14}$ & $\Z/ (2)$ & $2$ & $0$ & $6$ \\  
\hline  \hline  212  & $A\sb{1}+A\sb{3}+D\sb{14}$ & $\Z/ (2)$ & $2$ & $0$ & $4$ \\  
\hline  \hline  213  & $A\sb{4}+D\sb{14}$ & $(0)$ & $4$ & $2$ & $6$ \\  
\hline  \hline  214  & $A\sb{1}+A\sb{2}+D\sb{15}$ & $(0)$ & $4$ & $0$ & $6$ \\  
\hline  \hline  215  & $2\,A\sb{1}+D\sb{16}$ & $\Z/ (2)$ & $2$ & $0$ & $2$ \\  
\hline  \hline  216  & $A\sb{2}+D\sb{16}$ & $\Z/ (2)$ & $2$ & $1$ & $2$ \\  
\hline  \hline  217  & $A\sb{1}+D\sb{17}$ & $(0)$ & $2$ & $0$ & $4$ \\  
\hline  \hline  218  & $D\sb{18}$ & $(0)$ & $2$ & $0$ & $2$ \\  
\hline  \hline  219  & $3\,E\sb{6}$ & $\Z/ (3)$ & $2$ & $1$ & $2$ \\  
\hline  \hline  220  & $2\,A\sb{3}+2\,E\sb{6}$ & $(0)$ & $12$ & $0$ & $12$ \\  
\hline  \hline  221  & $A\sb{1}+A\sb{3}+2\,A\sb{4}+E\sb{6}$ & $(0)$ & $20$ & $0$ & $30$ \\  
\hline  \hline  222  & $A\sb{1}+A\sb{5}+2\,E\sb{6}$ & $\Z/ (3)$ & $2$ & $0$ & $6$ \\  
\hline  \hline  223  & $A\sb{2}+2\,A\sb{5}+E\sb{6}$ & $\Z/ (3)$ & $6$ & $0$ & $6$ \\  
\hline  \end{EllTab}  \vfill  \eject  
\begin{EllTab} 
\hline  \hline  224  & $2\,A\sb{2}+A\sb{3}+A\sb{5}+E\sb{6}$ & $\Z/ (3)$ & $6$ & $0$ & $12$ \\  
\hline  \hline  225  & $A\sb{3}+A\sb{4}+A\sb{5}+E\sb{6}$ & $(0)$ & $12$ & $0$ & $30$ \\  
\hline  \hline  226  & $A\sb{6}+2\,E\sb{6}$ & $(0)$ & $6$ & $3$ & $12$ \\  
\hline  \hline  227  & $A\sb{1}+A\sb{2}+A\sb{3}+A\sb{6}+E\sb{6}$ & $(0)$ & $6$ & $0$ & $84$ \\  
\cline{4-6}  &  &  & $12$ & $0$ & $42$ \\ 
\hline  \hline  228  & $2\,A\sb{1}+A\sb{4}+A\sb{ 6}+E\sb{6}$ & $(0)$ & $20$ & $10$ & $26$ \\  
\hline  \hline  229  & $A\sb{2}+A\sb{4}+A\sb{ 6}+E\sb{6}$ & $(0)$ & $18$ & $3$ & $18$ \\  
\hline  \hline  230  & $A\sb{1}+A\sb{5}+A\sb{6} +E\sb{6}$ & $(0)$ & $6$ & $0$ & $42$ \\  
\hline  \hline  231  & $A\sb{1}+A\sb{4}+A\sb{7}+E\sb{6}$ & $(0)$ & $2$ & $0$ & $120$ \\  
\hline  \hline  232  & $A\sb{5}+A\sb{7}+E\sb{6 }$ & $(0)$ & $6$ & $0$ & $24$ \\  
\hline  \hline  233  & $2\,A\sb{2}+A\sb{8}+E\sb{6}$ & $\Z/ (3)$ & $6$ & $3$ & $6$ \\  
\hline  \hline  234  & $2\,A\sb{1}+A\sb{2}+A\sb{8 }+E\sb{6}$ & $\Z/ (3)$ & $2$ & $0$ & $18$ \\  
\hline  \hline  235  & $A\sb{1}+A\sb{3}+A\sb{8}+E\sb{6}$ & $(0)$ & $12$ & $0$ & $18$ \\  
\hline  \hline  236  & $A\sb{ 4}+A\sb{8}+E\sb{6}$ & $(0)$ & $12$ & $3$ & $12$ \\  
\hline  \hline  237  & $A\sb{1}+A\sb{2}+A\sb{9} +E\sb{6}$ & $(0)$ & $12$ & $6$ & $18$ \\  
\hline  \hline  238  & $A\sb{3}+A\sb{9}+E\sb{6}$ & $(0)$ & $10$ & $0$ & $12$ \\  
\hline  \hline  239  & $2\,A\sb{1}+A\sb{10}+E\sb{6}$ & $(0)$ & $2$ & $0$ & $66$ \\  
\hline  \hline  240  & $A\sb{2}+A\sb{10}+E\sb{6}$ & $(0)$ & $6$ & $3$ & $18$ \\  
\hline  \hline  241  & $A\sb{1}+A\sb{11}+E\sb{6}$ & $(0)$ & $6$ & $0$ & $12$ \\  
\cline{3-6}  &  & $\Z/ (3)$ & $2$ & $0$ & $4$ \\  
\hline  \hline  242  & $A\sb{12}+E\sb{6}$ & $(0)$ & $4$ & $1$ & $10$ \\  
\hline  \hline  243  & $A\sb{3}+A\sb{4}+D\sb{5}+E\sb{6}$ & $(0)$ & $12$ & $0$ & $20$ \\  
\hline  \hline  244  & $A\sb{1}+A\sb{6}+D\sb{5}+E\sb{6}$ & $(0)$ & $2$ & $0$ & $84$ \\  
\hline  \hline  245  & $A\sb{7}+D\sb{5}+E\sb{6}$ & $(0)$ & $8$ & $0$ & $12$ \\  
\hline  \hline  246  & $D\sb{6}+2\,E\sb{6}$ & $(0)$ & $6$ & $0$ & $6$ \\  
\hline  \hline  247  & $A\sb{2}+A\sb{4}+D\sb{6}+E\sb{6}$ & $(0)$ & $6$ & $0$ & $30$ \\  
\hline  \hline  248  & $A\sb{6}+D\sb{6}+E\sb{6}$ & $(0)$ & $4$ & $2$ & $22$ \\  
\hline  \hline  249  & $A\sb{1}+A\sb{4}+D\sb{7}+E\sb{6}$ & $(0)$ & $4$ & $0$ & $30$ \\  
\hline  \hline  250  & $D\sb{5}+D\sb{7}+E\sb{6}$ & $(0)$ & $4$ & $0$ & $12$ \\  
\hline  \hline  251  & $A\sb{4}+D\sb{8}+E\sb{6}$ & $(0)$ & $8$ & $2$ & $8$ \\  
\hline  \hline  252  & $A\sb{1}+A\sb{2}+D\sb{9}+E\sb{6}$ & $(0)$ & $6$ & $0$ & $12$ \\  
\hline  \hline  253  & $A\sb{3}+D\sb{9}+E\sb{6}$ & $(0)$ & $4$ & $0$ & $12$ \\  
\hline  \hline  254  & $A\sb{1}+D\sb{11}+E\sb{6}$ & $(0)$ & $2$ & $0$ & $12$ \\  
\hline  \hline  255  & $D\sb{12}+E\sb{6}$ & $(0)$ & $4$ & $2$ & $4$ \\  
\hline  \hline  256  & $2\,A\sb{2}+2\,E\sb{7}$ & $(0)$ & $6$ & $0$ & $6$ \\  
\hline  \hline  257  & $A\sb{1}+A\sb{3}+2\,E\sb{7}$ & $\Z/ (2)$ & $2$ & $0$ & $4$ \\  
\hline  \end{EllTab}  \vfill  \eject  
\begin{EllTab} 
\hline  \hline  258  & $A\sb{4}+2\,E\sb{7}$ & $(0)$ & $4$ & $2$ & $6$ \\  
\hline  \hline  259  & $A\sb{1}+2\,A\sb{3}+A\sb{4}+E\sb{7}$ & $\Z/ (2)$ & $4$ & $0$ & $20$ \\  
\hline  \hline  260  & $2\,A\sb{2}+A\sb{3}+A\sb{4}+E\sb{7}$ & $(0)$ & $12$ & $0$ & $30$ \\  
\hline  \hline  261  & $2\,A\sb{3}+A\sb{5}+E\sb{7}$ & $\Z/ (2)$ & $4$ & $0$ & $12$ \\  
\hline  \hline  262  & $A\sb{1}+A\sb{2}+A\sb{3}+A\sb{5}+E\sb{7}$ & $\Z  /   (2)$ & $6$ & $0$ & $12$ \\  
\hline  \hline  263  & $2\,A\sb{1}+A\sb{4}+A\sb{5}+E\sb{7}$ & $\Z/ (2)$ & $8$ & $2$ & $8$ \\  
\hline  \hline  264  & $A\sb{2}+A\sb{4}+A\sb{5}+E\sb{7}$ & $(0)$ & $6$ & $0$ & $30$ \\  
\hline  \hline  265  & $A\sb{1}+2\,A\sb{2}+A\sb{6}+E\sb{7}$ & $(0)$ & $6$ & $0$ & $42$ \\  
\hline  \hline  266  & $A\sb{2}+A\sb{3}+A\sb{6}+E\sb{7}$ & $(0)$ & $4$ & $0$ & $42$ \\  
\hline  \hline  267  & $A\sb{1}+A\sb{4}+A\sb{6}+E\sb{7}$ & $(0)$ & $2$ & $0$ & $70$ \\  
\cline{4-6}  &  &  & $8$ & $2$ & $18$ \\  
\hline  \hline  268  & $A\sb{5}+A\sb{6}+E\sb{7}$ & $(0)$ & $4$ & $2$ & $22$ \\  
\hline  \hline  269  & $2\,A\sb{2}+A\sb{7}+E\sb{7}$ & $(0)$ & $6$ & $0$ & $24$ \\  
\hline  \hline  270  & $2\,A\sb{1}+A\sb{2}+A\sb{7}+E\sb{7}$ & $\Z/ (2)$ & $2$ & $0$ & $24$ \\  
\hline  \hline  271  & $A\sb{1}+A\sb{3}+A\sb{7}+E\sb{7}$ & $\Z/ (2)$ & $4$ & $0$ & $8$ \\  
\hline  \hline  272  & $A\sb{4}+A\sb{7}+E\sb{7}$ & $(0)$ & $6$ & $2$ & $14$ \\  
\hline  \hline  273  & $A\sb{1}+A\sb{2}+A\sb{8}+E\sb{7}$ & $(0)$ & $6$ & $0$ & $18$ \\  
\hline  \hline  274  & $A\sb{3}+A\sb{8}+E\sb{7}$ & $(0)$ & $4$ & $0$ & $18$ \\  
\hline  \hline  275  & $2\,A\sb{1}+A\sb{9}+E\sb{7}$ & $\Z/ (2)$ & $2$ & $0$ & $10$ \\  
\hline  \hline  276  & $A\sb{2}+A\sb{9}+E\sb{7}$ & $(0)$ & $6$ & $0$ & $10$ \\  
\cline{3-6}  &  & $\Z/ (2)$ & $4$ & $1$ & $4$ \\  
\hline  \hline  277  & $A\sb{1}+A\sb{10}+E\sb{7}$ & $(0)$ & $2$ & $0$ & $22$ \\  
\cline{4-6}  &  & & $6$ & $2$ & $8$ \\  
\hline  \hline  278  & $A\sb{11}+E\sb{7}$ & $(0)$ & $4$ & $0$ & $6$ \\  
\hline  \hline  279  & $D\sb{4}+2\,E\sb{7}$ & $\Z/ (2)$ & $2$ & $0$ & $2$ \\  
\hline  \hline  280  & $A\sb{2}+A\sb{4}+D\sb{5}+ E\sb{7}$ & $(0)$ & $6$ & $0$ & $20$ \\  
\hline  \hline  281  & $A\sb{1}+A\sb{5}+D\sb{5}+E\sb{7 }$ & $\Z/ (2)$ & $2$ & $0$ & $12$ \\  
\hline  \hline  282  & $A\sb{6}+D\sb{5}+E\sb{7}$ & $(0)$ & $6$ & $2$ & $10$ \\  
\hline  \hline  283  & $A\sb{2}+A\sb{3}+D\sb{6}+E\sb{7}$ & $\Z/ (2)$ & $4$ & $0$ & $6$ \\  
\hline  \hline  284  & $A\sb{5}+D\sb{6}+E\sb{7}$ & $\Z/ (2)$ & $4$ & $2$ & $4$ \\  
\hline  \hline  285  & $D\sb{5}+D\sb{6}+E\sb{7}$ & $\Z/ (2)$ & $2$ & $0$ & $4$ \\  
\hline  \hline  286  & $A\sb{1}+A\sb{3}+D\sb{7}+E\sb{7}$ & $\Z/ (2)$ & $4$ & $0$ & $4$ \\  
\hline  \hline  287  & $A\sb{4}+D\sb{7}+E\sb{7}$ & $(0)$ & $2$ & $0$ & $20$ \\  
\hline  \hline  288  & $A\sb{1}+A\sb{2}+D\sb{8}+E\sb{7}$ & $\Z/ (2)$ & $2$ & $0$ & $6$ \\  
\hline  \hline  289  & $A\sb{2}+D\sb{9}+E\sb{7}$ & $(0)$ & $4$ & $0$ & $6$ \\  
\hline  \hline  290  & $A\sb{1}+D\sb{10}+E\sb{7}$ & $\Z/ (2)$ & $2$ & $0$ & $2$ \\  
\hline  \end{EllTab} 
\vfill  \eject  
\begin{EllTab} 
\hline  \hline  291  & $D\sb{11}+E\sb{7}$ & $(0)$ & $2$ & $0$ & $4$ \\  
\hline  \hline  292  & $A\sb{2}+A\sb{3}+E\sb{6}+E\sb{7}$ & $(0)$ & $6$ & $0$ & $12$ \\  
\hline  \hline  293  & $A\sb{1}+A\sb{4}+E\sb{6}+E\sb{7}$ & $(0)$ & $2$ & $0$ & $30$ \\  
\hline  \hline  294  & $A\sb{5}+E\sb{6}+E\sb{7}$& $(0)$ & $6$ & $0$ & $6$ \\  
\hline  \hline  295  & $D\sb{5}+ E\sb{6}+E\sb{7}$ & $(0)$ & $2$ & $0$ & $12$ \\ 
\hline  \hline  296  & $2\,A\sb{1}+2\,E\sb{8}$ & $(0)$ & $2$ & $0$ & $2$ \\  
\hline  \hline  297  & $A\sb{2}+2\,E\sb{8}$ & $(0)$ & $2$ & $1$ & $2$ \\  
\hline  \hline  298  & $2\,A\sb{2}+2\,A\sb{3}+E\sb{8}$ & $(0)$ & $12$ & $0$ & $12$ \\  
\hline  \hline  299  & $2\,A\sb{1}+2\,A\sb{4}+E\sb{8}$ & $(0)$ & $10$ & $0$ & $10$ \\  
\hline  \hline  300  & $A\sb{1}+A\sb{2}+A\sb{3}+A\sb{4}+E\sb{8}$ & $(0)$ & $6$ & $0$ & $20$ \\  
\hline  \hline  301  & $2\,A\sb{5}+E\sb{8}$ & $(0)$ & $6$ & $0$ & $6$ \\  
\hline  \hline  302  & $A\sb{2}+A\sb{3}+A\sb{5}+E\sb{8}$ & $(0)$ & $6$ & $0$ & $12$ \\  
\hline  \hline  303  & $A\sb{1}+A\sb{4}+A\sb{5}+E\sb{8}$ & $(0)$ & $2$ & $0$ & $30$ \\  
\hline  \hline  304  & $2\,A\sb{2}+A\sb{6}+E\sb{8}$ & $(0)$ & $6$ & $3$ & $12$ \\  
\hline  \hline  305  & $2\,A\sb{1}+A\sb{2}+A\sb{6}+E\sb{8}$ & $(0)$ & $2$ & $0$ & $42$ \\  
\hline  \hline  306  & $A\sb{1}+A\sb{3}+A\sb{6}+E\sb{8}$ & $(0)$ & $6$ & $2$ & $10$ \\  
\hline  \hline  307  & $A\sb{4}+A\sb{6}+E\sb{8}$ & $(0)$ & $2$ & $1$ & $18$ \\  
\hline  \hline  308  & $A\sb{1}+A\sb{2}+A\sb{7}+E\sb{8}$ & $(0)$ & $2$ & $0$ & $24$ \\  
\hline  \hline  309  & $2\,A\sb{1}+A\sb{8}+E\sb{8}$ & $(0)$ & $2$ & $0$ & $18$ \\  
\hline  \hline  310  & $A\sb{2}+A\sb{8}+E\sb{8}$ & $(0)$ & $6$ & $3$ & $6$ \\  
\hline  \hline  311  & $A\sb{1}+A\sb{9}+E\sb{8}$ & $(0)$ & $2$ & $0$ & $10$ \\  
\hline  \hline  312  & $A\sb{10}+E\sb{8}$ & $(0)$ & $2$ & $1$ & $6$ \\  
\hline  \hline  313  & $2\,D\sb{5}+E\sb{8}$ & $(0)$ & $4$ & $0$ & $4$ \\  
\hline  \hline  314  & $A\sb{1}+A\sb{4}+D\sb{5}+E\sb{8}$ & $(0)$ & $2$ & $0$ & $20$ \\  
\hline  \hline  315  & $A\sb{5}+D\sb{5}+E\sb{8}$ & $(0)$ & $2$ & $0$ & $12$ \\  
\hline  \hline  316  & $2\,A\sb{ 2}+D\sb{6}+E\sb{8}$ & $(0)$ & $6$ & $0$ & $6$ \\  
\hline  \hline  317  & $A\sb{4}+D\sb{6}+E\sb{8}$ & $(0)$ & $4$ & $2$ & $6$ \\  
\hline  \hline  318  & $A\sb{1}+A\sb{2}+D\sb{7}+E\sb{8}$ & $(0)$ & $4$ & $0$ & $6$ \\  
\hline  \hline  319  & $A\sb{1}+D\sb{9}+E\sb{8}$ & $(0)$ & $2$ & $0$ & $4$ \\  
\hline  \hline  320  & $D\sb{10}+E\sb{8}$ & $(0)$ & $2$ & $0$ & $2$ \\  
\hline  \hline  321  & $A\sb{1}+A\sb{3}+E\sb{6}+E\sb{8}$ & $(0)$ & $2$ & $0$ & $12$ \\  
\hline  \hline  322  & $A\sb{4}+E\sb{6}+E\sb{8}$ & $(0)$ & $2$ & $1$ & $8$ \\  
\hline  \hline  323  & $D\sb{4}+E\sb{6}+E\sb{8}$ & $(0)$ & $4$ & $2$ & $4$ \\  
\hline  \hline  324  & $A\sb{1}+A\sb{2}+E\sb{7}+E\sb{8}$ & $(0)$ & $2$ & $0$ & $6$ \\  
\hline  \hline  325  & $A\sb{3}+E\sb{7}+E\sb{8}$ & $(0)$ & $2$ & $0$ & $4$ \\  
\hline  \end{EllTab}

\vfill\eject

\bibliographystyle{amsplain}

\begin{thebibliography}{10}
%
\bibitem {ATZ}
E.~Artal-Bartolo, H.~Tokunaga and D.~Q.~Zhang.
\textit{Miranda-Persson's problem on extremal elliptic K3 surfaces.}
preprint.
http://xxx.lanl.gov/list/math.AG, 9809065.
%
\bibitem{Bourbaki}
N.~Bourbaki.
{\'El\'ements de math\'ematique.
Groupes et alg\`ebres de Lie. Chapitre  IV-VI.}
Hermann, Paris, 1968.
%
\bibitem{CS}
J.~H.~Conway and  N.~J.~A.~Sloane. \textit{Sphere packings, lattices and groups.}
 Second edition.
 Grundlehren der Mathematischen Wissenschaften, \textbf{290},  Springer, New York,
1993.
%
\bibitem{Fujiki}
A.~Fujiki.  \textit{Finite automorphism groups of complex tori of
dimension two.} Publ.\ Res.\ Inst.\ Math.\ Sci. \textbf{24} (1988), no.~1, 1--97.
%
\bibitem {Kondo1}
S.~Kond\=o. 
\textit{Automorphisms of algebraic $K3$ surfaces which act trivially on Picard
groups.}
 J.\ Math.\ Soc.\ Japan \textbf{44} (1992), no.~1,
75--98. 
%
\bibitem{Kondo2}
S.~Kond\=o. \textit{Niemeier lattices, Mathieu groups, and
finite groups of symplectic automorphisms of $K3$ surfaces.} With an appendix by
Shigeru Mukai. Duke Math.\ J. \textbf{92} (1998), no. 3, 593--603. 
%
\bibitem {MP}
R.~Miranda and  U.~Persson. 
\textit{Mordell-Weil groups of extremal
elliptic $K3$ surfaces.}
 Problems in the theory of surfaces and their classification
(Cortona, 1988), Sympos.\ Math., XXXII, Academic Press, London, 1991, pp.~167--192.
%
%
\bibitem {Morrison}
D.~R.~Morrison. \textit{On $K3$ surfaces with large Picard number.}
Invent.\ Math. \textbf{75} (1984), no.~1, 105--121. 
%
\bibitem {Mukai}
S.~Mukai. \textit{Finite groups of automorphisms of $K3$
surfaces and the Mathieu group.} Invent.\ Math. \textbf{94} (1988), no.~1, 183--221. 
%
\bibitem {Nikulin1} 
V.~V.~Nikulin.
\textit{
Finite automorphism groups of K\"ahler $K3$ surfaces.}
Trans. Moscow Math. Soc. (1980), Issue 2,
pp.~71--135.
%
\bibitem {Nikulin2}
V.~V.~Nikulin. \textit{Integer symmetric bilinear forms and some of
their applications.}
Math. USSR Izvestija
\textbf{14} (1980), no.~1, 103--167.
%
\bibitem {Nishiyama} 
K.~Nishiyama. \textit{The Jacobian fibrations on some $K3$
surfaces and their Mordell-Weil groups.} Japan.\ J.\ Math. (N.S.) \textbf{22} (1996), no.~2,
293--347. 
%
\bibitem {Nori}
M.~V.~Nori. \textit{Zariski's conjecture and related problems.} Ann.\ Sci.\ \'Ecole Norm.\ Sup. (4) 
\textbf{16} (1983), no.~2, 305--344. 
%
\bibitem {PS}
I.~Piateskii-Shapiro and  I.~R.~Shafarevich.
\textit{A Torelli theorem for algebraic surfaces of type $K3$.}
Math.\ USSR Izv. \textbf{35} (1971), 530--572. 
%
\bibitem {Serre}
J.-P.~Serre. \textit{A course in arithmetic.} Graduate
Texts in Mathematics, \textbf{7}, Springer, New York, 1973.
%
\bibitem {Shioda-Inose}
T.~Shioda and  H.~Inose. \textit{On singular $K3$ surfaces.} Complex
analysis and algebraic geometry, Iwanami Shoten, Tokyo, 1977,
 pp.~119--136.
%
\bibitem {Todorov}
A.~N.~Todorov.
\textit{ Applications of the K\"ahler-Einstein-Calabi-Yau metric to
moduli of $K3$ surfaces}. 
Invent.\ Math. \textbf{61} (1980), no.~3, 251--265. 
%
\bibitem{Xiao}
G.~Xiao. \textit{ Galois covers between $K3$ surfaces.}
 Ann.\ Inst.\ Fourier (Grenoble) \textbf{46} (1996), no.~1, 73--88. 
%
\bibitem{Ye}
Q.~Ye.
\textit{On extremal elliptic $K3$ surfaces.}
preprint.
http://xxx.lanl.gov/abs/math.AG, 9901081 
%
\end{thebibliography}

\end{document}